\documentclass[12pt,a4paper]{article}
\usepackage[T1]{fontenc}
\usepackage{latexsym}
\usepackage{amsmath}
\usepackage{amsthm}
\usepackage{amssymb,enumerate}
\usepackage[usenames]{color}
\usepackage{graphicx}
\usepackage{fancybox}
\usepackage[normalem]{ulem}
\usepackage{hyperref}

\usepackage{color}

\newcommand{\N}{\mathbb{N}}

\newcommand{\R}{\mathbb{R}}
\newcommand{\C}{\mathbb{C}}
\newcommand{\HH}{\mathbb{H}}
\newcommand{\SMA}{\mathcal{S}^{\{\A\}}_{\{\M\}}}

\newcommand{\CM}{C_{\{\M\}}(0,\infty)}

\newcommand{\M}{{\boldsymbol{M}}}
\newcommand{\m}{{\boldsymbol{m}}}
\newcommand{\NN}{{\boldsymbol{N}}}
\newcommand{\A}{{\boldsymbol{A}}}
\newcommand{\hM}{\widehat{{\boldsymbol{M}}}}
\newcommand{\hNN}{\widehat{{\boldsymbol{N}}}}
\newcommand{\hA}{\widehat{{\boldsymbol{A}}}}

\newcommand{\id}{\operatorname{id}}

\newcommand{\dx}{{\rm d}x }

\newcommand{\supp}{\operatorname{supp}}
\newcommand{\lc}{\operatorname{(lc)}}
\newcommand{\dc}{\operatorname{(dc)}}
\newcommand{\mg}{\operatorname{(mg)}}
\newcommand{\alg}{\operatorname{(alg)}}
\newcommand{\nq}{\operatorname{(nq)}}
\newcommand{\snq}{\operatorname{(snq)}}
\newcommand{\sm}{\operatorname{(sm)}}

\newcommand{\om}{\omega}
\newcommand{\ga}{\gamma}

\theoremstyle{plain}
\newtheorem{theorem}{Theorem}[section]
\newtheorem{proposition}[theorem]{Proposition}
\newtheorem{corollary}[theorem]{Corollary}
\newtheorem{lemma}[theorem]{Lemma}
\newtheorem{prop}[theorem]{Proposition}

\newtheorem*{coro*}{Corollary}
\theoremstyle{definition}
\newtheorem{remark}[theorem]{Remark}
\newtheorem{definition}[theorem]{Definition}

\newtheorem{examples}[theorem]{Examples}
\numberwithin{equation}{section}

\makeatletter
\DeclareRobustCommand\widecheck[1]{{\mathpalette\@widecheck{#1}}}
\def\@widecheck#1#2{%
    \setbox\z@\hbox{\m@th$#1#2$}%
    \setbox\tw@\hbox{\m@th$#1%
       \widehat{%
          \vrule\@width\z@\@height\ht\z@
          \vrule\@height\z@\@width\wd\z@}$}%
    \dp\tw@-\ht\z@
    \@tempdima\ht\z@ \advance\@tempdima2\ht\tw@ \divide\@tempdima\thr@@
    \setbox\tw@\hbox{%
       \raise\@tempdima\hbox{\scalebox{1}[-1]{\lower\@tempdima\box
\tw@}}}%
    {\ooalign{\box\tw@ \cr \box\z@}}}
\makeatother

\begin{document}

\title{The Stieltjes moment problem in Gelfand-Shilov spaces defined by weight sequences in the absence of derivation closedness}

\author{Javier Jim\'enez-Garrido \and Ignacio Miguel-Cantero \and Javier Sanz \and Gerhard Schindl}
\date{\today}

\maketitle

\begin{abstract}

The Stieltjes moment problem is studied in a new framework within the general Gelfand-Shilov spaces defined via weight sequences. The novelty consists of allowing for a naturally larger target space for the moment mapping, which sends a function to its sequence of Stieltjes moments. The motivation comes from a recent version of the Borel-Ritt theorem, concerning the surjectivity of the Borel mapping in Carleman-Roumieu ultraholomorphic classes in sectors, whose defining weight sequence is subject to the condition, weaker than derivation closedness, of having shifted moments. The injectivity and surjectivity of the moment mapping in this new setting is studied and, in some cases, characterized. Finally, results are provided for general weight sequences of fast and regular enough growth when the condition of shifted moments fails to hold.
\par\medskip

\noindent Key words: Stieltjes moment problem, Gelfand-Shilov spaces, growth and regularity conditions for sequences.\par
\medskip
\noindent 2020 MSC: 30E05, 44A60, 46F05.
\end{abstract}

\section{Introduction}\label{intro}

This paper is a new contribution to the study of the Stieltjes moment problem in the context of Gelfand-Shilov spaces of Roumieu type defined by weight sequences, first introduced in~\cite{G-S} and which we define now.
$\mathcal{S}(0,\infty)$ will denote the subspace of the Schwartz class of rapidly decreasing complex-valued smooth functions in $\R$ with support in $[0,\infty)$, and $\N_0=\{0,1,2,\dots\}$. Given two sequences of positive real numbers $\M=(M_p)_{p\in\N_0}$ and $\A=(A_p)_{p\in\N_0}$, we consider the Gelfand-Shilov spaces of Roumieu type
$\mathcal{S}_{\{\M\}}^{\{\A\}}(0,\infty)$ and $\mathcal{S}_{\{\M\}}(0,\infty)$,
consisting of all $\varphi \in \mathcal{S}(0,\infty)$  such that there exists $h>0$ with
$$
\sup_{p,q \in \N_0} \sup_{x \in \R} \frac{|x^p\varphi^{(q)}(x)|}{h^{p+q}M_pA_q} < \infty
$$
and
$$
\sup_{p\in \N_0} \sup_{x \in \R} \frac{|x^p\varphi^{(q)}(x)|}{h^{p}M_p} < \infty \qquad \text{for every }q \in \N_0,
$$
respectively. It is clear that $\mathcal{S}_{\{\M\}}^{\{\A\}}(0,\infty)\subset \mathcal{S}_{\{\M\}}(0,\infty)$, and that for every $\varphi \in \mathcal{S}_{\{\M\}}(0,\infty)$ the sequence of Stieltjes moments $(\mu_p(\varphi))_{p\in\N_0} = (\int_0^\infty x^p \varphi(x) \dx)_{p \in \N_0}$ is well defined and has a restricted growth.
In case $\M$ is derivation closed, that is, $M_{p+1} \leq C_0H^{p+1}M_p$ for every  $p \in \N_0$ and some $C_0>0$ and $H \geq 1$, it is easy to check that the moment sequence belongs to
$
\Lambda_{\{\M\}} = \{ (c_p)_{p \in \N_0}\colon \sup_{p \in \N_0} \frac{|c_p|}{h^p M_p} < \infty \text{ for some }h>0 \}.
$
The standard Stieltjes problem in this context consists then in the study of the surjectivity and injectivity of the Stieltjes moment mapping $\mathcal{M}\colon\varphi\mapsto (\mu_p(\varphi))_{p\in\N_0}$, when defined on either $\mathcal{S}_{\{\M\}}^{\{\A\}}(0,\infty)$ or $\mathcal{S}_{\{\M\}}(0,\infty)$ and with range $\Lambda_{\{\M\}}$.
Surjectivity has been studied in a series of papers, always resting on ideas of A.~L.~Dur\'an and R.~Estrada~\cite{D-E} that combine the Fourier transform with Borel-Ritt-like theorems from asymptotic analysis, see S.-Y.~Chung, D.~Kim and Y.~Yeom~\cite[Thm.\ 3.1]{ChungKimYeom} for $\M=(p!^{\alpha})_{p\in\N_0}$ (the Gevrey sequences) whenever $\alpha>2$, and A. Lastra and the third author~\cite{L-S08,L-S09} for $\mathcal{S}_{\{(p!M_p)_p\}}(0,\infty)$ and general strongly regular sequences $\M$ whose growth index $\gamma(\M)$ is strictly greater than 1 (see Section \ref{sect-prelim-mod} for the definitions).
Subsequently, A. Debrouwere and the first and third  authors~\cite{DebrouwereetalMoment} improved and completed these results by including the spaces $\mathcal{S}_{\{(p!M_p)_p\}}^{\{(p!A_p)_p\}}(0,\infty)$ in their considerations, by dropping some hypotheses on $\M$ (specially moderate growth, stronger than derivation closedness), and by also studying the injectivity of the Stieltjes moment mapping.
The new key tools were a better understanding of the meaning of the different growth conditions usually imposed on the sequence $\M$ and their expression in terms of indices of O-regular variation, as developed in~\cite{JimenezSanzSchindlIndex}, and the enhanced information obtained in~\cite{JG-S-S} about the injectivity and surjectivity of the asymptotic Borel mapping on Carleman-Roumieu ultraholomorphic classes in sectors defined by sequences $\M$ subject to minimal conditions. Finally, A. Debrouwere~\cite{momentsdebrouwere} completely characterized the surjectivity and the existence of global right inverses for the moment mapping in Gelfand-Shilov spaces of both Roumieu and Beurling type under derivation closedness. His technique does not rest on Borel-Ritt-like theorems, but relates the problem to the surjectivity and existence of global right inverses for the Borel mapping in Carleman ultradifferentiable classes, already characterized by H.-J. Petzsche~\cite{Petzsche88}.

The main aim of this paper is the study of the Stieltjes moment problem in a new framework, allowing for a naturally larger target space for the moment mapping. The motivation comes from a recent version of the Borel-Ritt theorem, concerning the surjectivity of the Borel mapping in Carleman-Roumieu ultraholomorphic classes in sectors, whose defining weight sequence $\M$ is subject to a condition $\sm$ much weaker than derivation closedness. This constructive version uses a truncated Laplace transform with an integral kernel $e(z)$ given in terms of a so-called optimal flat function within the ultraholomorphic class. The key fact is that $\sm$ characterizes the equivalence of $\M_{+1}:=(M_{p+1})_p$ and the Stieltjes moment sequence of $e$ restricted to the positive real line, what makes the procedure work. So, it becomes natural (see Propositions~\ref{prop.dc.charact.LambdaM} and~\ref{prop.sm.charact.LambdaM+1}) to change the target space into the larger one $\Lambda_{\{\M_{+1}\}}$, and study again the injectivity and surjectivity in this new setting. As the technique in~\cite{momentsdebrouwere} does not seem to apply, we have recovered the technique in~\cite{L-S09}, resting on the construction of local right inverses for the moment mapping. This requires a careful study of the action of the Fourier transform under this new condition $\sm$ (Proposition~\ref{Fourier-contention-dm}), and the adaptation of some auxiliary results which were already useful in previous frameworks.

In a final section we explore the situation when condition $\sm$ fails to hold and, so, the target space can no longer be $\Lambda_{\{\M_{+1}\}}$. It turns out that the same accurate estimates leading us to the formulation of condition $\sm$ are the key for proposing the new target space, associated with a new sequence $\widetilde{\M}$ introduced in~\eqref{eq.defWidetildeM}. In this framework, injectivity and surjectivity results can be proved in a similar way as long as the weight sequence $\M$ has a fast and regular growth guaranteeing that $\widetilde{\M}$ is again a weight sequence.

The paper is organized as follows. In the preliminary Section \ref{sect-prelim-mod} we gather the main facts needed regarding sequences, ultraholomorphic classes, the asymptotic Borel mapping $\mathcal{B}$ and growth indices related to the injectivity and surjectivity of $\mathcal{B}$, Gelfand-Shilov spaces and the Fourier and Laplace transform.

Section~\ref{sect-Stiel-mom-problem-dm} contains the main results. Firstly, the use of optimal flat functions in ultraholomorphic classes, whose existence was shown in~\cite{JimenezMiguelSanzSchindlOptFlat}, allows to determine the appropriate target space in the moment problem according to whether $\dc$ or $\sm$ is satisfied. After some auxiliary results,  Theorem~\ref{inj-stiel-dm} characterizes the injectivity of the Stieltjes moment mapping under condition $\sm$ for $\M$, while  Theorem~\ref{surj-stiel-modified-dm} studies the surjectivity problem and its connection to the existence of local right inverses for $\mathcal{M}$ with a uniform scaling of the parameter defining the Banach spaces under consideration.
Finally, Section~\ref{sect-Stiel-mom-problem-no-cond} consists of the study of the problem in the absence of $\sm$.

\section{Preliminaries}\label{sect-prelim-mod}

We set $\N=\{1,2,\dots\}$, $\N_0=\{0,1,2,\dots\}$.
$\M = (M_p)_{p \in \N_0}$  will stand for a sequence of positive real numbers with $M_0 = 1$. The {\it sequence of quotients} of $\M$ is  $\m=(m_p)_{p\in\N_0}$ given by
$m_p:=\frac{M_{p+1}}{M_p}$, $p\in \N_0$. The knowledge of $\M$ amounts to that of $\m$, since $M_p=m_0\cdots m_{p-1}$, $p\in\N$. We will denote by small letters the quotients of a sequence given by the corresponding capital letters.
We set $\widehat{\M}:=(p!M_p)_{p \in \N_0}$ and $\widecheck{\M}:=(M_p/p!)_{p \in \N_0}$; $(\widehat{m}_p)_p$ and $(\widecheck{m}_p)_p$ will denote the corresponding sequences of quotients: $\widehat{m}_p=(p+1)m_p$, $\widecheck{m}_p=m_p/(p+1)$ for $p\in\N_0$.

We shall use the following conditions on  sequences $\M$:

\begin{itemize}
\item[$\lc$] $\M$ is \emph{log-convex} if $M^2_p \leq M_{p-1}M_{p+1}$, $p \in \N$. Equivalently, $(m_p)_{p\in\N_0}$ is nondecreasing.
\item[$\sm$] $\M$ has \emph{shifted moments} if $\log(m_{p+1}/m_p) \leq C_0H^{p+1}$, $p \in \N_0$, for some $C_0>0$ and $H>1$.
\item[$\dc$] $\M$ is \emph{derivation closed} if $M_{p+1} \leq C_0H^{p+1}M_p$, $p \in \N_0$, for some $C_0>0$ and $H>1$.
\item[$\mg$] $\M$ has \emph{moderate growth} if $M_{p+q} \leq C_0H^{p+q}M_pM_q$, $p,q \in \N_0$, for some $C_0>0$ and $H>1$.
\item[$\alg$] $\M$ satisfies the \emph{algebrability condition} if $M_pM_q\leq C_1^{p+q}M_{p+q} $, $p,q \in \N_0$, for some $C_1\ge 1$.
\item[$\nq$] $\M$ is \emph{non-quasianalytic} if $\displaystyle \sum_{p=0}^\infty \frac{1}{(p+1)m_p} < \infty$.
\item[$\snq$] $\M$ is \emph{strongly non-quasianalytic} if $\displaystyle \sum_{q=p}^\infty \frac{1}{(q+1)m_q} \leq  \frac{C}{m_p}$, $p \in \N_0$, for some $C > 0$.
\end{itemize}

In the classical work of H.~Komatsu~\cite{Komatsu}, the properties $\lc$, $\dc$ and $\mg$ are denoted by $(M.1)$, $(M.2)'$ and $(M.2)$, respectively, while $\nq$ and $\snq$ for $\M$ are the same as properties $(M.3)'$ and $(M.3)$ for $\widehat{\M}$, respectively.

If $\M$ is $\lc$, not only $(m_p)_{p\in\N_0}$ but also $(M_p^{1/p})_{p\in\N}$ is nondecreasing, and $(M_p)^{1/p}\leq m_{p-1}$ for every $p\in\N$; moreover, $\lim_{p\to\infty} (M_p)^{1/p}= \infty$ if and only if $\lim_{p\to\infty} m_p= \infty$.
The sequence $\M$ is said to be a \emph{weight sequence} if it satisfies $\lc$ and $m_p \nearrow \infty$ as $p \to \infty$.

\begin{remark}\label{rem_propertiesMandstability}
All the previously listed properties are preserved when passing from $\M$ to $\widehat{\M}$. However, only $\sm$, $\dc$, $\mg$ and $\alg$ are generally kept when going from $\M$ to $\widecheck{\M}$.
It is also clear that $\snq$ implies $\nq$, and that $\mg$ implies $\dc$. Moreover, as soon as $a_0:=\inf_{p\in\N_0}m_p>0$ (and, in particular, if $\M$ is $\lc$) one can check that $\dc$ implies $\sm$. The statements concerning $\sm$ have been proved in~\cite[Sect. 2.1]{JimenezMiguelSanzSchindlSurjectSM}.
\end{remark}

For a sequence $\M$ satisfying $\nq$ and such that $\hM$ satisfies $\lc$, one easily proves that $m_p \rightarrow \infty$ as $p \to \infty$.
According to E.~M.~Dyn'kin~\cite{Dynkin80}, if $\M$ is a weight sequence and satisfies (dc), we say $\hM$ is \emph{regular}. Following V.~Thilliez~\cite{Thilliez03}, if $\M$ satisfies $\lc$, $\mg$ and $\snq$, we say $\M$ is \emph{strongly regular}; in this case $\M$ is a weight sequence, and the corresponding $\hM$ is regular.

\begin{examples}
\begin{itemize}
\item[(i)] The Gevrey sequences $(p!^{\alpha})_{p}$ ($\alpha>0$) are strongly regular, and their perturbations $(p!^{\alpha}\prod_{j=0}^p\log^{\beta}(e+j))_{p}$ ($\alpha>0$, $\beta\in\R$) also are after suitably modifying a finite number of their terms, if necessary. They appear all through the study of formal power series solutions to differential and difference equations.
\item[(ii)] The sequences $\M=(q^{p^\alpha})_{p}$ ($q>1$, $0<\alpha\le 2$) and $\M=(p^{\tau p^\sigma})_{p}$ ($\tau>0$, $1<\sigma<2$) are such that $\hM$ is regular, but they are not strongly regular since they do not satisfy $\mg$. In particular, the $q$-Gevrey sequences $\M=(q^{p^2})_{p}$ ($q>1$) appear in the study of $q$-difference equations.
\item[(iii)] The weight sequences $\M=(q^{p^\alpha})_{p}$ ($q>1$, $\alpha>2$) and $\M=(p^{\tau p^\sigma})_{p}$ ($\tau>0$, $\sigma\ge 2$) do not satisfy $\dc$, so that $\hM$ is not regular, but they still satisfy $\sm$.
\item[(iv)] The rapidly growing weight sequence $\M=(q^{p^p})_{p}$ ($q>1$; $M_0:=1$) does not satisfy $\sm$. However, the Stieltjes moment problem can still be treated in this case by a modified argument included in Section~\ref{sect-Stiel-mom-problem-no-cond}.
\end{itemize}

\end{examples}

For later use, we recall that $\lc$ (together with $M_0=1$) implies that
\begin{equation}\label{eq-conseq-lc-prod}
M_pM_q\le M_{p+q},\quad p,q\in\N_0,
\end{equation}
and so the condition $\alg$ is satisfied.

Following Komatsu~\cite{Komatsu}, the relation $\M \subset \NN$ means that there are $C,h > 0$ such that $M_p \leq Ch^pN_p$ for all $p \in \N_0$. We say $\M$ and $\NN$ are \emph{equivalent}, denoted  $\M \approx \NN$, if $\M \subset \NN$ and $\NN \subset \M$. It is straightforward that $\sm$ (see~\cite[Lemma 2.4]{JimenezMiguelSanzSchindlSurjectSM}), $\dc$, $\mg$ and $\alg$ are stable under equivalence for general sequences, and the same can be deduced for $\nq$ and $\snq$ for weight sequences by indirect methods, as these two last conditions characterize the non injectivity and the surjectivity, respectively, of the Borel map in Carleman ultradifferentiable classes, by the classical Denjoy-Carleman theorem (see, for example,~\cite{RudinRCA}) and the results of H.-J. Petzsche~\cite[Cor. 3.2]{Petzsche88} (see~\cite[Cor. 3.14]{JimenezSanzSchindlIndex} for a direct proof of a more general statement about the stability of $\snq$).

The \emph{associated weight function} of $\M$ is defined as
$$
\omega_{\M}(t) := \sup_{p \in \N_0} \log \frac{t^p}{M_p}, \qquad t > 0,
$$
and $\omega_{\M}(0) := 0$. Another associated function is
\begin{equation*}
	h_\M(t):=\inf_{p\in\N_0}M_p t^p, \quad t>0;\ h_\M(0):=0.
\end{equation*}
The functions $h_\M$ and $\omega_\M$ are related by
\begin{equation}\label{functionhequ2}
	h_\M(t)=\exp(-\omega_\M(1/t)),\quad t>0.
\end{equation}
It is well known that $\om_{\M}(t)<\infty$ for every $t>0$ if, and only if, $\lim_{p\to\infty}M_p^{1/p}=\infty$, which will be our standard assumption when $\om_{\M}$ enters our considerations. In particular, this condition is satisfied for weight sequences. In~\cite[Prop.~3.2]{Komatsu} we find that, for a weight sequence $\M$,
\begin{equation}\label{eq.MpfromomegaM}
	M_p=\sup_{t>0}t^p\exp(-\om_{\M}(t))=\sup_{t>0}t^p h_{\M}(1/t),\quad p\in\N_0,
\end{equation}
and it is plain to check that
\begin{equation}\label{eq.expression_h_M}
h_\M(t)=M_{p+1}t^{p+1} \quad\textrm{for }\frac{1}{m_{p+1}}\le t <\frac{1}{m_p},\ p\in\N_0;\quad h_\M(t)=1\quad\textrm{for }t\ge\frac{1}{m_0},
\end{equation}
so that an explicit expression is also available for $\om_{\M}$.
Moreover, $h_\M$ and $\om_{\M}$ are nondecreasing and continuous.

\subsection{Ultraholomorphic classes on the upper half-plane and the asymptotic Borel mapping}\label{subsect-ultraholom-dm}

Given an open set $\Omega$ in the complex plane $\C$, we denote by
$\mathcal{O}(\Omega)$ the space of holomorphic functions in $\Omega$. We write $\HH$ for the open upper half-plane of~$\C$.

Let $\M=(M_p)_p$ be a sequence of positive real numbers. In the definition of the following classes of functions or sequences it is not relevant that $M_0=1$, but, as said before, this will be assumed unless stated otherwise. For $h>0$
we define $\mathcal{A}_{\M,h}(\HH)$ as the space consisting of all $f \in \mathcal{O}(\HH)$ such that
$$
\|f\|_{\M,h}:=\sup_{p \in \N_0} \sup_{z \in \HH} \frac{|f^{(p)}(z)|}{h^pM_p} < \infty.
$$
$(\mathcal{A}_{\M,h}(\HH),\|\cdot\|_{\M,h})$ is a Banach space, and we set $\mathcal{A}_{\{\M\}}(\HH) = \bigcup_{h>0} \mathcal{A}_{\M,h}(\HH)$, endowed with the corresponding $(LB)$ space structure.

The next result follows from the fact that the elements of $\mathcal{A}_{\M,h}(\HH)$ together with all their derivatives are Lipschitz on $\HH$.

\begin{lemma}\label{extension-dm}
Let $\M$ be a sequence and let $f \in \mathcal{A}_{\M,h}(\HH)$ for some $h > 0$. Then,
\begin{equation*}
f_p(x) = \lim_{z \in \HH, z \to x} f^{(p)}(z)  \in \C
\end{equation*}
exists for all $x \in \R$ and $p \in \N_0$. Moreover, $f_0 \in C^\infty(\R)$, $f_0^{(p)} = f_p$ for all $p \in \N$, and
$$
\sup_{p \in \N_0} \sup_{x \in \R} \frac{|f_0^{(p)}(x)|}{h^pM_p} < \infty.
$$
\end{lemma}

In the sequel, we shall frequently write $f(x) = \lim_{z \in \HH, z \to x} f(z)$ for $x \in \R$, and so $f$ and all its derivatives can be considered to be continuous on $\overline{\HH}$ and satisfy the same global estimates there.

Let $\M$ be a sequence. For $h>0$ we define $\Lambda_{\M,h}$  as the space consisting of all sequences $(c_p)_p \in \C^{\N_0}$ such that
$$
|(c_p)_p|_{\M,h}:=\sup_{p \in \N_0} \frac{|c_p|}{h^pM_p} < \infty.
$$
$(\Lambda_{\M,h},|\cdot|_{\M,h})$ is a Banach space, and  $\Lambda_{\{\M\}}:= \bigcup_{h > 0} \Lambda_{\M,h}$ is the corresponding $(LB)$ space. The \emph{asymptotic Borel mapping} is defined as
$$
\mathcal{B}: \mathcal{A}_{\{\M\}}(\HH) \rightarrow \Lambda_{\{\M\}},\ \  f \mapsto (f^{(p)}(0))_{p},
$$
which is well-defined by Lemma \ref{extension-dm}. Since $|\mathcal{B}(f)|_{\M,h}\le \|f\|_{\M,h}$ for every $f\in\mathcal{A}_{\M,h}(\HH)$, it is linear and continuous both between $\mathcal{A}_{\M,h}(\HH)$ and $\Lambda_{\M,h}$ for every $h>0$, and between $\mathcal{A}_{\{\M\}}(\HH)$ and $\Lambda_{\{\M\}}$.
An updated account on the injectivity and surjectivity of the asymptotic Borel mapping on various ultraholomorphic classes defined on arbitrary sectors may be found in the works of the authors \cite{JG-S-S,JimenezSanzSchindlSurjectDC,JimenezMiguelSanzSchindlOptFlat,JimenezMiguelSanzSchindlSurjectSM}. There, two indices $\ga(\M)$ and $\om(\M)$, associated to the sequence $\M$, play a prominent role.
In~\cite[Ch.~2]{PhDJimenez} and \cite[Sect.~3]{JimenezSanzSchindlIndex}, the connections between these indices, the growth properties usually imposed on  sequences, and the theory of O-regular variation, have been thoroughly studied. We summarize here some facts.

Recall that a sequence $(c_p)_{p}$ is \emph{almost increasing} if there exists $a>0$ such that $c_p\leq a c_q $ for all $ q\geq p$. The first index, introduced by V. Thilliez~\cite[Sect.\ 1.3]{Thilliez03} for strongly regular sequences, may be defined for any weight sequence $\M$ as
$$\ga(\M):=\sup\{\mu>0\, | \,(m_{p}/(p+1)^\mu)_{p}\hbox{ is almost increasing} \}\in[0,\infty].
$$
For  $\beta>0$ we say that $\M$  satisfies $(\ga_{\beta})$ if there exists $C>0$ such that
\begin{equation*}
\sum^\infty_{q=p} \frac{1}{(m_q)^{1/\beta}}\leq \frac{C (p+1) }{(m_p)^{1/\beta}},  \qquad p\in\N_0.
\end{equation*}
Then one has that
$$ \ga(\M)=\sup\{\beta>0  \, | \, \M \,\, \text{satisfies} \,\, (\ga_{\beta})  \}.$$

For future reference we gather some more information.
\begin{lemma}\label{lemma.PropertiesGammaIndex}
Let $\M$ be a weight sequence and $\beta>0$. The following statements hold:
\begin{enumerate}
 \item[(i)] $\ga(\M)>0$ if and only if $\M$ satisfies $\snq$.
 \item[(ii)] $\ga(\M)>\beta$ if and only if $\M$ satisfies $(\ga_\beta)$.
 \item[(iii)] $\ga(\M)>\beta$ implies that $(p!^{\beta})_p\subset\M$.
\end{enumerate}
\end{lemma}

\begin{remark}\label{rem.weight.almost.incr}
Every weight sequence $\M$ is almost increasing: since eventually $m_p\ge 1$, we have that $M_p\le M_{p+1}$ if $p\ge p_0$ for suitable $p_0$. Then $M_p\le a M_q$ for any $p\le q$ if we choose $a:=\max_{p\le q\le p_0}M_p/M_q\ge 1$.
\end{remark}

The surjectivity of the asymptotic Borel mapping in a half-plane can be characterized as follows. Note that the condition $\gamma(\M)>1$ amounts, in view of~(ii) and the easy equality $\gamma(\hM)=\gamma(\M)+1$, to the fact that $\hM$ satisfies $(\gamma_2)$, which is the condition appearing in~\cite[Thm. 7.4.(b)]{momentsdebrouwere}.

\begin{theorem}[\cite{momentsdebrouwere}, Theorem 7.4.(b)]\label{teor.Andreas}
Let $\M$ be a weight sequence satisfying $\dc$. The following are equivalent:
\begin{enumerate}[(i)]
\item $\mathcal{B}\colon \mathcal{A}_{\{\M\}}(\HH)\to\Lambda_{\{\M\}}$ is surjective.
\item $\gamma(\M)>2$.
\end{enumerate}
\end{theorem}

The second index $\om(\M)$ is given by
$$\om(\M):= \displaystyle\liminf_{p\to\infty} \frac{\log(m_{p})}{\log(p)}\in[0,\infty],
$$
and it turns out that
\begin{align*}
\om(\M)&=\sup\{\mu>0\, | \, \sum^\infty_{p=0} \frac{1}{(m_p)^{1/\mu}}<\infty \}\\
&=\sup\{\mu>0 \, | \, \sum^\infty_{p=0} \frac{1}{((p+1)m_p)^{1/(\mu+1)}} < \infty \}.
\end{align*}

Concerning the injectivity of the asymptotic Borel mapping, we have the next result.
\begin{theorem} \emph{(\cite[Thm.\ 12]{Salinas}, \cite[Thm.\ 3.4]{JG-S-S})}\label{inj-borel-dm} Let $\M$ be a weight sequence. Then, $\mathcal{B}: \mathcal{A}_{\{\M\}}(\HH) \rightarrow \Lambda_{\{\M\}}$ is injective if and only if
$$
\sum_{p = 0}^\infty \frac{1}{m_p^{1/2}} = \infty,
$$
which in turn implies that $\om(\M)\le 2$.
\end{theorem}

Finally, we mention that if $\M$ is a weight sequence, the inequality $\gamma(\M)\le\omega(\M)$ holds, and so the asymptotic Borel mapping  $\mathcal{B}: \mathcal{A}_{\{\M\}}(\HH) \rightarrow \Lambda_{\{\M\}}$ is never bijective, see also~\cite[Thm.\ 3.17]{JG-S-S}.

\subsection{Gelfand-Shilov spaces and the Fourier transform}\label{GSspaceswithdm}
Let $\M$ and $\A$ be sequences of positive real numbers. For $h>0$ we define $\mathcal{S}_{\M,h}^{\A,h}(\R)$ as the space consisting of all $\varphi \in C^\infty(\R)$ such that
$$
s_{\M,h}^{\A,h}(\varphi):=\sup_{p,q \in \N_0} \sup_{x \in \R} \frac{|x^p\varphi^{(q)}(x)|}{h^{p+q}M_pA_q} < \infty.
$$
$(\mathcal{S}_{\M,h}^{\A,h}(\R),s_{\M,h}^{\A,h})$ is a Banach space. If $M_p^{1/p}\to\infty$ as $p\to\infty$, and in particular if $\M$ is a weight sequence, notice that $\varphi \in C^\infty(\R)$ belongs to  $\mathcal{S}_{\M,h}^{\A,h}(\R)$ if and only if
$$
\sup_{q \in \N_0} \sup_{x \in \R} \frac{|\varphi^{(q)}(x)|e^{\om_{\M}(|x|/h)}}{h^{q}A_q} < \infty
$$
or, in other words, for every $q \in \N_0$ one has
\begin{equation}\label{eq.fibelongsSMAviahM}
|\varphi^{(q)}(x)|\le s_{\M,h}^{\A,h}(\varphi) h^{q}A_qe^{-\om_{\M}(|x|/h)}= s_{\M,h}^{\A,h}(\varphi)h^{q}A_qh_{\M}(h/|x|),\ \ x\in\R
\end{equation}
(with the agreement that $h_{\M}(h/|x|)=1$ for $x=0$, according to the fact that $\om_{\M}(0)=0$).

We set $\SMA(\R) = \bigcup_{h>0} \mathcal{S}_{\M,h}^{\A,h}(\R)$, which is an $(LB)$ space.

Analogously, we define  $\mathcal{S}_{\M,h}(\R)$, $h> 0$, as the space consisting of all $\varphi \in C^\infty(\R)$ such that, for all $q\in\N_0$,
$$
s_{\M,h}^{q}(\varphi):=\sup_{p\in \N_0} \sup_{x \in \R} \frac{|x^p\varphi^{(q)}(x)|}{h^{p}M_p} < \infty.
$$
$s_{\M,h}^{q}$ is a seminorm, $(\mathcal{S}_{\M,h}(\R),(s_{\M,h}^{q})_{q\in\N_0})$ is a Fr\'echet space, and we set $\mathcal{S}_{\{\M\}}(\R) = \bigcup_{h>0} \mathcal{S}_{\M,h}(\R)$, endowed with its natural $(LF)$ space structure. We also define
$$
\SMA(0,\infty) := \{ \varphi \in \SMA(\R) \, | \, \supp \varphi \subseteq [0,\infty)\}
$$
and
$$
\mathcal{S}_{\{\M\}}(0,\infty) := \{ \varphi \in \mathcal{S}_{\{\M\}}(\R) \, | \, \supp \varphi \subseteq [0,\infty)\},
$$
whose relative topologies from their ambient spaces coincide, as long as $\M$ is a weight sequence, with the corresponding $(LB)$ and $(LF)$ structures obtained from the similarly defined Banach subspaces $\mathcal{S}_{\M,h}^{\A,h}(0,\infty)$ or Fr\'echet subspaces $\mathcal{S}_{\M,h}(0,\infty)$, see~\cite[Lemma 3.3 and page 24]{momentsdebrouwere}.

Observe that $\SMA(\R)\subset \mathcal{S}_{\{\M\}}(\R)$ and $\SMA(0,\infty)\subset \mathcal{S}_{\{\M\}}(0,\infty)$. If $\A$ satisfies $\lc$, then $\SMA(0,\infty)$ is non-trivial (i.e., it contains non identically zero functions) if and only if $\sum_{p=0}^\infty 1/a_p<\infty$, as follows from the Denjoy-Carleman theorem.

In the remainder of this subsection we investigate the image of the spaces $\SMA(\R)$ and $\SMA(0,\infty)$
under the Fourier transform (cf.\ \cite[Sect.\ IV.6]{G-S}), which we define as follows:
$$
\mathcal{F}(\varphi)(\xi) = \widehat{\varphi}(\xi) = \int_{-\infty}^\infty \varphi(x) e^{ix\xi} \dx, \qquad \varphi \in L^1(\R).
$$
With this definition it is well-known that
\begin{equation}\label{eq.InverseFourierTransform}
\mathcal{F}^{-1}(\varphi)(\xi)=\frac{1}{2\pi}\mathcal{F}(\varphi)(-\xi),
\ \ \varphi\in\mathcal{S}(\R),\ \xi\in\R,
\end{equation}
where $\mathcal{S}(\R)$ is the Schwartz space of complex-valued, rapidly decreasing smooth functions.

In our next statements, the sequence
$$\M_{+1}:=(M_{p+1})_{p\in\N_0}$$
will play a prominent role. Note that its first term $M_1$ will not be generally equal to 1.

\begin{proposition}\label{Fourier-contention-dm}
Let $\M$ be a weight sequence satisfying $\sm$, and $\A$ be either an almost increasing sequence, or a sequence such that $\liminf_{p\to\infty}A_p^{1/p}>0$ and $\hA$ satisfies $\alg$. Then:
\begin{enumerate}[(i)]
\item There exists $a>0$ such that for every $h\ge 1$ one has $\mathcal{F}(\mathcal{S}^{\hA,h}_{\M,h}(\R))\subset \mathcal{S}^{\M_{+1},ah}_{\hA,ah}(\R)$, and $\mathcal{F}\colon \mathcal{S}^{\hA,h}_{\M,h}(\R)\to \mathcal{S}^{\M_{+1},ah}_{\hA,ah}(\R)$ is continuous.
\item  $\mathcal{F}(\mathcal{S}^{\{\hA\}}_{\{\M\}}(\R))\subset \mathcal{S}^{\{\M_{+1}\}}_{\{\hA\}}(\R)$ and $\mathcal{F}\colon\mathcal{S}^{\{\hA\}}_{\{\M\}}(\R)\to \mathcal{S}^{\{\M_{+1}\}}_{\{\hA\}}(\R)$ is continuous.
\item The two previous statements are valid when replacing $\mathcal{F}$ by $\mathcal{F}^{-1}$.
\end{enumerate}
\end{proposition}
\begin{proof}
The statements about $\mathcal{F}^{-1}$ will be valid if the ones about $\mathcal{F}$ are, due to~\eqref{eq.InverseFourierTransform}.
(ii) is immediate from (i), which we prove now. Let $h \geq 1$ and $\varphi \in \mathcal{S}_{\M,h}^{\hA,h}(\R)$ be arbitrary. Then,
\begin{equation}\label{eq.SMhA}
\sup_{x \in \R} |x^p\varphi^{(q)}(x)| \leq s_{\M,h}^{\hA,h}(\varphi)h^{p+q}M_pq!A_q, \qquad p,q \in \N_0.
\end{equation}
After applying Leibniz's theorem for parametric integrals and integration by parts, we get
\begin{equation}\label{eq.estimate-Fourier-first}
\sup_{\xi \in \R} |\xi^q\widehat{\varphi}^{(p)}(\xi)| \leq \sum_{j=0}^{\min\{p,q\}} \binom{q}{j} \frac{p!}{(p-j)!} \int_{-\infty}^\infty |x^{p-j} \varphi^{(q-j)}(x)| \dx.
\end{equation}
We split each of the integrals into five intervals, and use~\eqref{eq.SMhA} in order to estimate the first of them:
\begin{align}
\int_{-\infty}^{-hm_{p+1-j}} |x^{p+2-j} \varphi^{(q-j)}(x)|\frac{\dx}{x^2}&\le
s_{\M,h}^{\hA,h}(\varphi)h^{p+q+2-2j}M_{p+2-j}(q-j)!A_{q-j} \frac{1}{hm_{p+1-j}}\nonumber\\
&=s_{\M,h}^{\hA,h}(\varphi)h^{p+q+1-2j}M_{p+1-j}(q-j)!A_{q-j}.\label{eq.estimate-Fourier-extreme-intervals}
\end{align}
One can proceed similarly for the integral over $(hm_{p+1-j},\infty)$. Regarding the interval $(hm_{p-j},hm_{p+1-j})$, one uses~\eqref{eq.fibelongsSMAviahM}, then~\eqref{eq.expression_h_M} and finally $\sm$ to obtain
\begin{align}
\int_{hm_{p-j}}^{hm_{p+1-j}} x^{p-j} |\varphi^{(q-j)}(x)|\dx &\le
s_{\M,h}^{\hA,h}(\varphi)h^{q-j}(q-j)!A_{q-j} \int_{hm_{p-j}}^{hm_{p+1-j}}x^{p-j} h_{\M}(h/x)\dx\nonumber\\
&=
s_{\M,h}^{\hA,h}(\varphi)h^{q-j}(q-j)!A_{q-j} \int_{hm_{p-j}}^{hm_{p+1-j}}x^{p-j} \frac{h^{p+1-j}M_{p+1-j}}{x^{p+1-j}}\dx\nonumber\\
&=
s_{\M,h}^{\hA,h}(\varphi)h^{p+q+1-2j}M_{p+1-j}(q-j)!A_{q-j} \log(\frac{m_{p+1-j}}{m_{p-j}})\nonumber\\
&\le C_0s_{\M,h}^{\hA,h}(\varphi)H^{p+1-j}h^{p+q+1-2j}M_{p+1-j}(q-j)!A_{q-j}.
\label{eq.estimate-Fourier-middle-intervals}
\end{align}
The integral on the interval $(-hm_{p+1-j},-hm_{p-j})$ is treated analogously. Finally, again~\eqref{eq.SMhA} provides
\begin{align}
\int_{-hm_{p-j}}^{hm_{p-j}} |x^{p-j}\varphi^{(q-j)}(x)|\dx &\le
2hm_{p-j}s_{\M,h}^{\hA,h}(\varphi)h^{p+q-2j}M_{p-j}(q-j)!A_{q-j}\nonumber\\
&=
2s_{\M,h}^{\hA,h}(\varphi)h^{p+q+1-2j}M_{p+1-j}(q-j)!A_{q-j}.
\label{eq.estimate-Fourier-central-interval}
\end{align}
Taking~\eqref{eq.estimate-Fourier-extreme-intervals}, \eqref{eq.estimate-Fourier-middle-intervals} and \eqref{eq.estimate-Fourier-central-interval} into~\eqref{eq.estimate-Fourier-first}, we obtain
\begin{equation}\label{eq.estimatesFourierTransfGlobal}
\sup_{\xi \in \R} |\xi^q\widehat{\varphi}^{(p)}(\xi)|\leq \sum_{j=0}^{\min\{p,q\}} \binom{q}{j}\binom{p}{j} j! 2s_{\M,h}^{\hA,h}(\varphi)h^{p+q+1-2j}M_{p+1-j}(q-j)!A_{q-j} (2+C_0H^{p+1-j}).
\end{equation}
$\M$ is almost increasing (see Remark~\ref{rem.weight.almost.incr}). In case $\A$ also is, there exists $D\ge 1$ such that $A_j\le DA_p$ and $M_j\le DM_p$ for all $j\le p$. This fact allows us to write
\begin{align*}
\sup_{\xi \in \R} |\xi^q\widehat{\varphi}^{(p)}(\xi)|&\leq
2s_{\M,h}^{\hA,h}(\varphi)D^2h^{p+q+1}M_{p+1}q!A_{q}(2+C_0H)
\sum_{j=0}^{p}\binom{p}{j}h^{-2j}H^{p-j}\\
&= 2(2+C_0H)D^2hs_{\M,h}^{\hA,h}(\varphi)\left(\frac{1}{h^2}+H\right)^p h^{p+q}M_{p+1}q!A_{q}\\
&\le 2(2+C_0H)D^2hs_{\M,h}^{\hA,h}(\varphi)\left((1+H)h\right)^{p+q} M_{p+1}q!A_{q}
\end{align*}
for all $p,q \in \N_0$.

In case $\hA$ satisfies $\alg$, we have that $(q-j)!A_{q-j}\le C_1^qq!A_q/(j!A_j)$ for $0\le j\le q$ and some $C_1\ge 1$. If moreover $\liminf_{p\to\infty}A_p^{1/p}>0$, there exists $c>0$ such that $A_p\ge c^p$ for every $p$. Going back to~\eqref{eq.estimatesFourierTransfGlobal} and using that  $\binom{q}{j}\le 2^q$ for every $0\le j\le q$, we get
\begin{align*}
\sup_{\xi \in \R} |\xi^q\widehat{\varphi}^{(p)}(\xi)|&\leq
2s_{\M,h}^{\hA,h}(\varphi)Dh^{p+q+1}2^qC_1^q M_{p+1}q!A_{q}(2+C_0H)
\sum_{j=0}^{p}\binom{p}{j}h^{-2j}c^{-j}H^{p-j}\\
&= 2(2+C_0H)Dhs_{\M,h}^{\hA,h}(\varphi)\left(\frac{1}{ch^2}+H\right)^p2^q C_1^qh^{p+q}M_{p+1}q!A_{q}\\
&\le 2(2+C_0H)Dhs_{\M,h}^{\hA,h}(\varphi) \left(2C_1\Big(\frac{1}{c}+H\Big)h\right)^{p+q}M_{p+1}q!A_{q}
\end{align*}
for all $p,q \in \N_0$. Hence, we have proved the first statement with $a=H+1$ in the first case, and $a=2C_1(H+1/c)$ in the second one.
\end{proof}

\begin{remark}
In the previous result we can impose $\alg$ to be satisfied, equivalently, by either $\A$ or $\hA$. We have preferred the second option because, in some of its consequences (see, for example, Lemma~\ref{multiplication-dm} and Theorems~\ref{inj-stiel-dm} and~\ref{surj-stiel-modified-dm}), we will ask indeed $\hA$ to be a weight sequence, and so it immediately satisfies $\alg$.
\end{remark}

Resting on Proposition \ref{Fourier-contention-dm}, the next result can be shown in a similar way as the corresponding implication in \cite[Prop.\ 2.1]{C-C-K}.
\begin{proposition}\label{Fourier-image-dm} Let $\M$ be a weight sequence satisfying $\sm$, and $\A$ be either an almost increasing sequence, or a sequence such that $\liminf_{p\to\infty}A_p^{1/p}>0$ and $\hA$ satisfies $\alg$. If $\psi \in \mathcal{S}^{\{\M_{+1}\}}_{\{\hA\}}(\R)$ and there is $\Psi: \overline{\HH} \rightarrow \C$ satisfying the following conditions:
\begin{itemize}
\item[$(i)$] $\Psi_{|\R} = \psi$.
\item[$(ii)$] $\Psi$ is continuous on $\overline{\HH}$ and analytic on $\HH$.
\item[$(iii)$] $\lim_{\zeta \in \overline{\HH}, \zeta \to \infty} \Psi(\zeta) = 0$,
\end{itemize}
then $\psi \in \mathcal{F}(\mathcal{S}^{\{\hA\}}_{\{\M\}}(0,\infty))$.
\end{proposition}

\subsection{The Laplace transform}\label{sect-Laplace-dm} Let $\M$ be a sequence. We define $C_{\M,h}(0,\infty)$ as the space consisting of all $\varphi \in C((0,\infty))$ such that
$$
s_{\M,h}^0(\varphi)=\sup_{p \in \N_0} \sup_{x \in (0,\infty)} \frac{x^p|\varphi(x)|}{h^{p}M_p} < \infty.
$$
$(C_{\M,h}(0,\infty),s_{\M,h}^0)$ is a Banach space.
If $\M$ is a weight sequence and $\varphi\in C_{\M,h}(0,\infty)$,
similarly as in~\eqref{eq.fibelongsSMAviahM} we have
\begin{equation}\label{eq.fibelongsCMviahM}
|\varphi(x)|\le s_{\M,h}^0(\varphi)e^{-\om_{\M}(|x|/h)}= s_{\M,h}^0(\varphi)h_{\M}(h/|x|),\ \ x>0.
\end{equation}

We set $\CM = \bigcup_{h>0} C_{\M,h}(0,\infty)$, and endow it with its natural $(LB)$ space structure. Note that no condition is imposed on the support of the elements of $\CM$.

The \emph{Laplace transform} of $\varphi \in \CM$ is defined as
$$
\mathcal{L}(\varphi)(\zeta) =  \int_0^\infty \varphi(x)e^{ix\zeta}\dx, \qquad \zeta \in \overline{\HH}.
$$
\begin{remark}\label{rem.ContinuousInclusions}
Let $\M$ and $\A$ be sequences. We have $\SMA(0,\infty)\subset\mathcal{S}_{\{\M\}}(0,\infty)\subset\CM$ with continuous inclusions, since $\mathcal{S}_{\M,h}^{\A,h}(0,\infty)\subset\mathcal{S}_{\M,h}(0,\infty) \subset C_{\M,h}(0,\infty)$ for every $h>0$, the norm in $C_{\M,h}(0,\infty)$ enters the family of seminorms defining the topology of $\mathcal{S}_{\M,h}(0,\infty)$, and
$$
s_{\M,h}^q(\varphi)\le h^qA_qs_{\M,h}^{\A,h}(\varphi),\ \ \varphi\in\mathcal{S}_{\M,h}^{\A,h}(0,\infty),\ q\in\N_0.
$$
Note that $\mathcal{L}(\varphi)_{|\R} = \widehat{\varphi}$ for all $\varphi \in \SMA(0,\infty)$.
\end{remark}
\begin{lemma}\label{inj-lapl-dm}
Let $\M$ be a weight sequence satisfying $\sm$. Then, for every $h>0$ one has $\mathcal{L}(C_{\M,h}(0,\infty))\subset \mathcal{A}_{\M_{+1},Hh}(\HH)$, where $H>1$ is the constant appearing in $\sm$, and $\mathcal{L}\colon C_{\M,h}(0,\infty)\to \mathcal{A}_{\M_{+1},Hh}(\HH)$ is continuous. So, the mapping
$
\mathcal{L} : \CM \rightarrow \mathcal{A}_{\{\M_{+1}\}}(\HH)
$
is well-defined and continuous, and it is moreover injective.
\end{lemma}
\begin{proof}
Suppose $\varphi\in\CM$, and choose $h>0$ such that~\eqref{eq.fibelongsCMviahM} holds. Given
$\zeta\in\HH$ and $p\in\N_0$, since $\Re(ix\zeta)=-x\Im(\zeta)< 0$ for every $x> 0$, we have
\begin{align}
|(\mathcal{L}(\varphi))^{(p)}(\zeta)|&\le  \int_0^\infty x^p|\varphi(x)|\dx\label{eq.bound_Laplace_transform}\\
&\le \int_0^{hm_p}\! x^p |\varphi(x)|\dx+
s_{\M,h}^0(\varphi)\int_{hm_p}^{hm_{p+1}}\! x^p \frac{h^{p+1}M_{p+1}}{x^{p+1}}\dx
		+\int_{hm_{p+1}}^\infty \!\! \frac{x^{p+2}|\varphi(x)|}{x^2}\dx\nonumber\\
		&\le s_{\M,h}^0(\varphi)\left(hm_ph^pM_p+ h^{p+1}M_{p+1}\log\left(\frac{m_{p+1}}{m_p}\right)
		+h^{p+2}M_{p+2}\frac{1}{hm_{p+1}}\right)\nonumber\\
&\le s_{\M,h}^0(\varphi)h^{p+1}M_{p+1}\left(2+C_0H^{p+1}\right)\le
(2+C_0H)hs_{\M,h}^0(\varphi)(Hh)^pM_{p+1},\nonumber
\end{align}
where in the next-to-last inequality $\sm$ has been applied. Hence, $\mathcal{L}$ is well-defined and continuous from $C_{\M,h}(0,\infty)$ into $\mathcal{A}_{\M_{+1},hH}(\HH)$.

The proof of injectivity can be found in~\cite[Lemma 2.10]{DebrouwereetalMoment}.
\end{proof}

\begin{remark}\label{rem.Laplacewithdc}
If we suppose that $\M$ satisfies the stronger condition $\dc$ instead of $\sm$, then the Laplace transform $\mathcal{L}$ sends $\CM$ into $\mathcal{A}_{\{\M\}}(\HH)$. This is easily seen by splitting the integral in~\eqref{eq.bound_Laplace_transform} into only two subintervals, $(0,hm_{p})$ and $(hm_p,\infty)$ and estimating similarly.
\end{remark}

\section{A modified Stieltjes moment problem in Gelfand-Shilov spaces}\label{sect-Stiel-mom-problem-dm}
Let $\M$ be a weight sequence. The $p$-th moment, $p \in \N_0$, of an element $\varphi \in \CM$ is defined as
$$
\mu_p(\varphi) := \int_0^\infty x^p \varphi(x) \dx.
$$
The formula
\begin{equation}\label{eq.moments.deriv.Laplace}
\mathcal{L}(\varphi)^{(p)}(0) = i^p \mu_p(\varphi), \qquad \varphi \in \CM, p \in \N_0,
\end{equation}
guarantees, according to Lema~\ref{inj-lapl-dm}, that whenever $\M$ satisfies $\sm$ the \emph{Stieltjes moment mapping}
$$
\mathcal{M}: \CM \rightarrow \Lambda_{\{\M_{+1}\}};\ \varphi \mapsto (\mu_p(\varphi))_p
$$
is well-defined and continuous. Indeed, for every $h>0$ and $\varphi\in C_{\M,h}(0,\infty)$ one has
$$
|\mathcal{M}(\varphi)|_{\M_{+1},Hh}\le (2+C_0H)hs_{\M,h}^0(\varphi).
$$
However, if $\M$ satisfies $\dc$, Remark~\ref{rem.Laplacewithdc} shows that $\mathcal{M}$ sends $\CM$ into $\Lambda_{\{\M\}}$. This latter situation was studied in~\cite{DebrouwereetalMoment}, while the former one is our objective now. In order to stress the relevance of conditions $\dc$ and $\sm$ for the Stieltjes moment problem to be well-posed with one or another target space, we need to introduce the concept of optimal flat functions in ultraholomorphic classes in unbounded sectors of the Riemann surface of the logarithm $\mathcal{R}$. We refer to~\cite{JimenezMiguelSanzSchindlOptFlat} for further details.

\begin{definition}\label{optimalflatdef}
Let $\M$ be a weight sequence and $S\subset\mathcal{R}$ be an unbounded sector bisected by direction $d=0$, i.e., by the positive real line $(0,+\infty)$. A holomorphic function $G\colon S\to\C$ is called an \emph{optimal $\{\M\}$-flat function} in $S$ if:
	\begin{itemize}
		\item[$(i)$] There exist $K_1,K_2>0$ such that for all $x>0$,
		\begin{equation}\label{optimalflatleft}
			K_1h_{\M}(K_2x)\le G(x).
		\end{equation}
		
		\item[$(ii)$] There exist $K_3,K_4>0$ such that for all $z\in S$, one has
		\begin{equation}\label{optimalflatright}
			|G(z)|\le K_3h_{\M}(K_4|z|).
		\end{equation}
	\end{itemize}
\end{definition}

We note that the estimate in~\eqref{optimalflatright} amounts to the fact that
$$|G(z)|\le K_3K_4^pM_p|z|^p,\qquad p\in\N_0,\ z\in S,
$$
what exactly says that $G$ admits uniform $\{\M\}$-asymptotic expansion in $S$ given by the null power series, and this is the meaning of \emph{$\{\M\}$-flatness}. The inequality~\eqref{optimalflatleft} makes the function optimal in a sense, as its rate of decrease on the positive real axis when $t$ tends to 0 is accurately specified by the function $h_{\M}$.
The following result is interesting for our discussion.

\begin{proposition}(\cite[Prop. 3.10]{JimenezMiguelSanzSchindlOptFlat})\label{prop.ExistOptFlatFunct}
	Let $\M$ be a weight sequence with $\ga(\M)>0$. Then, for any $0<\gamma<\ga(\M)$ there exists an optimal $\{\M\}$-flat function in the sector $S_\ga=\{z\in\mathcal{R}\colon|\arg(z)|<\pi\ga/2\}$.
\end{proposition}

If $G$ is an optimal $\{\M\}$-flat function in $S_\ga$, we define the kernel function $e\colon S_\ga\to\C$ given by
\begin{equation*}
	e(z):=G\left(\frac{1}{z}\right),\quad z\in S_\ga.
\end{equation*}
Because of~\eqref{optimalflatleft} and~\eqref{optimalflatright}, we have that
\begin{equation}\label{eq.Bounds_e_sector}
K_1h_{\M}\left(\frac{K_2}{x}\right)\le e(x)\le K_3h_{\M}\left(\frac{K_4}{x}\right),\quad x>0,
\end{equation}
and according to~\eqref{eq.fibelongsCMviahM}, we see that (the restriction to $(0,\infty)$ of) $e$ belongs to $\CM$.

The following result, partially obtained in \cite[Prop. 3.11]{JimenezMiguelSanzSchindlOptFlat}, shows the key role of such kernel functions. We include the whole proof for the reader's convenience.

\begin{proposition}\label{prop.dc.charact.LambdaM}
Let $\M$ be a weight sequence with $\ga(\M)>0$. Then,
$\mathcal{M}(\CM)\subset\Lambda_{\{\M\}}$ if, and only if, $\M$ satisfies $\dc$.
\end{proposition}
\begin{proof}
As indicated in Remark~\ref{rem.Laplacewithdc}, the condition is sufficient. Conversely, suppose now that $\mathcal{M}(\CM)\subset\Lambda_{\{\M\}}$. Consider an optimal $\{\M\}$-flat function $G$ in a suitably narrow sector $S$ bisected by the positive real axis, and let $e$ be the corresponding kernel function. Since (the restriction to $(0,\infty)$ of) $e\in\CM$, there exists $C,h>0$ such that
$\mu_p(e)\le Ch^pM_p$ for every $p\in\N_0$. On the other hand, by the left-hand inequalities in~\eqref{eq.Bounds_e_sector} and the monotonicity of $h_{\M}$, for every $p\in\N_0$ and $s>0$ we may estimate
$$
\mu_p(e)\ge \int_0^s t^p e(t)\,dt \ge K_1 \int_0^s t^p h_{\M}\left(\frac{K_2}{t}\right)\,dt\ge
K_1 h_{\M}\left(\frac{K_2}{s}\right)\frac{s^{p+1}}{p+1}.
$$
Then, by~\eqref{eq.MpfromomegaM} we deduce that
$$
\mu_p(e)\ge \frac{K_1}{p+1} \sup_{s>0}h_{\M}\left(\frac{K_2}{s}\right)s^{p+1}=
\frac{K_1}{p+1} K_2^{p+1}M_{p+1}\ge
K_1K_2\left(\frac{K_2}{2}\right)^p M_{p+1}.
$$
From the estimates for $\mu_p(e)$ from above and below we deduce that $\dc$ is satisfied.
\end{proof}

Similarly, we have the following characterization, again partially included in~\cite{JimenezMiguelSanzSchindlSurjectSM}.

\begin{proposition}\label{prop.sm.charact.LambdaM+1}
Let $\M$ be a weight sequence with $\ga(\M)>0$. Then,
$\mathcal{M}(\CM)\subset\Lambda_{\{\M_{+1}\}}$ if, and only if, $\M$ satisfies $\sm$.
\end{proposition}

\begin{proof} As previously said, Lemma~\ref{inj-lapl-dm} implies that the condition is sufficient. Suppose now that $\mathcal{M}(\CM)\subset\Lambda_{\{\M_{+1}\}}$, and consider $e\in\CM$ as before, so that there exist $C,h>0$ such that
$\mu_p(e)\le Ch^pM_{p+1}$ for every $p\in\N_0$.
On the other hand, by the left-hand inequalities in~\eqref{eq.Bounds_e_sector}, for every $p\in\N_0$ we have that
	$$
	\mu_p(e)\ge \int_{K_2m_p}^{K_2m_{p+1}} t^p e(t)\,dt \ge K_1 \int_{K_2m_p}^{K_2m_{p+1}} t^p h_{\M}\left(\frac{K_2}{t}\right)\,dt =K_1K_2^{p+1}M_{p+1}\log\left(\frac{m_{p+1}}{m_p}\right),
	$$
where the last equality is a consequence of~\eqref{eq.expression_h_M}. Analogously as before, the estimates for $\mu_p(e)$ from above and below imply that $\sm$ is satisfied.
\end{proof}

\begin{remark}
Let $\M$ be a weight sequence with $\ga(\M)>0$. From the previous two results one easily concludes that $\M$ satisfies $\dc$ (respectively, $\sm$) if, and only if, for every kernel function $e$ associated with an optimal $\{\M\}$-flat function $G$ in some $S_\ga$ one has $(\mu_p(e))_p\approx \M$ (resp. $(\mu_p(e))_p\approx \M_{+1}$).
\end{remark}

We will reduce, via the Laplace transform, the study of the injectivity and surjectivity of the Stieltjes moment mapping in this new setting to their counterparts for the asymptotic Borel mapping (Theorems \ref{inj-borel-dm} and
\ref{teor.Andreas}), as it was already done by A. L. Dur\'an and R. Estrada in~\cite{D-E}, and later on by several authors \cite{C-C-K, L-S08,L-S09,DebrouwereetalMoment}.

The next lemma provides an auxiliary function, already appearing in the work \cite{D-E} and later adapted to our needs, see~\cite{DebrouwereetalMoment}. We set $\HH_{-1} = \{ z \in \C \, | \, \Im m \, z > -1 \}$.

\begin{lemma}(\cite[Lemma 3.1]{DebrouwereetalMoment})\label{aux-dm} Let $\A$ be a sequence satisfying $\nq$, and such that $\hA$ is a weight sequence. Then, there is $G_0 \in \mathcal{O}(\HH_{-1})$ satisfying the following conditions:
\begin{itemize}
\item[$(i)$] $G_0$ does not vanish on $\HH_{-1}$.
\item[$(ii)$] $\displaystyle \sup_{z \in \HH_{-1}} |G_0(z)|e^{\om_{\widehat{A}}(|z|)} < \infty$.
\item[$(iii)$] $\displaystyle \sup_{p \in \N}\sup_{x \in \R} \frac{|G_0^{(p)}(x)|e^{\om_{\widehat{A}}(|x|/2)}}{2^pp!} < \infty$.
\end{itemize}
\end{lemma}

Proposition \ref{Fourier-image-dm} and Lemma \ref{aux-dm} imply the following general result. A similar proof for classical Gelfand-Shilov spaces (i.e., when $\M$ is a Gevrey sequence $(p!^\alpha)_p$ with $\alpha>1$) can be found in~\cite[Prop. 4.13]{L-S08}, and a similar statement for strongly regular sequences $\M$ with $\ga(\M)>1$ was included in~\cite[Prop. 6.6]{L-S09} without proof. Another version, disregarding continuity and under stronger conditions than the ones imposed here, appeared in~\cite[Lemma 3.3]{DebrouwereetalMoment}. Note that, as indicated in Subsection~\ref{GSspaceswithdm}, the conditions imposed on the sequence $\A$ are precisely those guaranteeing that the class $\mathcal{S}^{\{\hA\}}_{\{\M\}}(0,\infty)$ is non-trivial.

\begin{lemma}\label{multiplication-dm} Let $\M$ be a sequence satisfying $\alg$ and such that $(p!)_p\subset \M$ (equivalently, $\liminf_{p\to\infty}\widecheck{M}^{1/p}_p>0$),
and let $\A$ be a sequence satisfying $\nq$ and such that $\hA$ is a weight sequence.
Consider the function $G_0$ from Lemma \ref{aux-dm}. Then, there exists $a>0$ such that for every $h\ge 1$ and for every $f \in \mathcal{A}_{\M_{+1},h}(\HH)$ one has $(fG_0)|_{\R} \in\mathcal{S}^{\M_{+1},ah}_{\hA,ah}(\R)$, and the map so defined is continuous from $\mathcal{A}_{\M_{+1},h}(\HH)$ into $\mathcal{S}^{\M_{+1},ah}_{\hA,ah}(\R)$.

\noindent If, additionally, $\M$ is a weight sequence  and satisfies also $\sm$, then for every $f \in \mathcal{A}_{\{\M_{+1}\}}(\HH)$ we have $(fG_0)|_{\R} \in \mathcal{F}(\mathcal{S}^{\{\hA\}}_{\{\M\}}(0,\infty))$.
\end{lemma}
\begin{proof}
According to Lemma~\ref{aux-dm}.$(iii)$, \eqref{functionhequ2} and the definition of $h_{\hA}$, there exists $C>0$ such that for every $p,j\in\N_0$ and $x\in\R\setminus\{0\}$ one has
\begin{equation}\label{eq.estimates.deriv.G_0}
|G_0^{(j)}(x)|\le C2^jj!h_{\hA}(2/|x|)\le C2^jj!\widehat{A}_p\left(\frac{2}{|x|}\right)^p,
\end{equation}
while $|G_0^{(j)}(0)|\le C2^jj!$ for every $j\in\N_0$.
Hence, for $h\ge 1$, $f \in \mathcal{A}_{\M_{+1},h}(\HH)$, $x\in\R$ and $p,q\in\N_0$ we get
\begin{align*}
|x^p(fG_0)^{(q)}(x)|&= \left|x^p\sum_{k=0}^q\binom{q}{k}f^{(k)}(x)G_0^{(q-k)}(x)\right|\\
&\le C\|f\|_{\M_{+1},h}2^p\widehat{A}_p \sum_{k=0}^q\binom{q}{k}h^kM_{k+1}2^{q-k}(q-k)!.
\end{align*}
Let $C_1$ be the constant in $\alg$ for $\M$. Also, by hypothesis there exists $B>0$ such that $M_p\ge B^pp!$ for every $p\in\N_0$, so
\begin{align*}
|x^p(fG_0)^{(q)}(x)|&\le
C\|f\|_{\M_{+1},h}2^p\widehat{A}_pC_1^{q+1}M_{q+1} \sum_{k=0}^q\binom{q}{k}h^k2^{q-k}\frac{(q-k)!}{M_{q-k}}\\
&\le C\|f\|_{\M_{+1},h}2^p\widehat{A}_pC_1^{q+1}M_{q+1} \sum_{k=0}^q\binom{q}{k}h^k\left(\frac{2}{B}\right)^{q-k}\\
&\le
CC_1\|f\|_{\M_{+1},h} \left(2C_1(1+2/B)h\right)^{p+q}\widehat{A}_pM_{q+1}.
\end{align*}
Hence, the first statement is proved with $a=2C_1(1+2/B)$. The second assertion stems directly from Proposition~\ref{Fourier-image-dm}.
\end{proof}

We are ready to study the injectivity and surjectivity of the Stieltjes moment mapping.

\begin{theorem}\label{inj-stiel-dm} Let $\M$ be a weight sequence satisfying $\sm$ and $(p!)_p\subset \M$, and let $\A$ be a sequence satisfying $\nq$ and such that $\hA$ is a weight sequence. Then, the following statements are equivalent:
\begin{itemize}
\item[$(i)$] $\displaystyle \sum_{p = 0}^\infty \frac{1}{m_p^{1/2}} = \infty$.
\item[$(ii)$] $\mathcal{B}: \mathcal{A}_{\{\M_{+1}\}}(\HH) \rightarrow \Lambda_{\{\M_{+1}\}}$ is injective.
\item[$(iii)$] $\mathcal{M}: \CM \rightarrow \Lambda_{\{\M_{+1}\}}$ is injective.
\item[$(iv)$] $\mathcal{M}: \mathcal{S}_{\{\M\}}(0,\infty) \rightarrow \Lambda_{\{\M_{+1}\}}$ is injective.
\item[$(v)$] $\mathcal{M}: \mathcal{S}^{\{\hA\}}_{\{\M\}}(0,\infty) \rightarrow \Lambda_{\{\M_{+1}\}}$ is injective.
\end{itemize}
\end{theorem}
\begin{proof}
$(i) \Rightarrow (ii)$: By Theorem \ref{inj-borel-dm}.

$(ii) \Rightarrow (iii)$: Let $\varphi \in \CM$ be such that $\mu_p(\varphi) = 0$ for all $p \in \N_0$. By Lemma \ref{inj-lapl-dm} we have that $\mathcal{L}(\varphi) \in \mathcal{A}_{\{\M_{+1}\}}(\HH)$.  Moreover,
$
\mathcal{L}(\varphi)^{(p)}(0) = i^p \mu_p(\varphi) = 0$ for all $p \in \N_0$ and, thus,
$\mathcal{L}(\varphi) \equiv 0$. Since $\mathcal{L}$ is injective (Lemma \ref{inj-lapl-dm}), we obtain that $\varphi \equiv 0$.

$(iii) \Rightarrow (iv) \Rightarrow (v)$: Obvious, see Remark~\ref{rem.ContinuousInclusions}.

$(v) \Rightarrow (i)$: In view of Theorem \ref{inj-borel-dm} it suffices to show that $\mathcal{B}: \mathcal{A}_{\{\M_{+1}\}}(\HH) \rightarrow \Lambda_{\{\M_{+1}\}}$ is injective. Let $f \in \mathcal{A}_{\{\M_{+1}\}}(\HH)$ be such that $f^{(p)}(0) = 0$ for all $p \in \N_0$. Consider the function $G_0$ from Lemma \ref{aux-dm}. By Lemma \ref{multiplication-dm} we have that $(fG_0)|_{\R} = \widehat{\varphi}$ for some $\varphi \in \mathcal{S}^{\{\hA\}}_{\{\M\}}(0,\infty)$. Observe that
$$
\mu_p(\varphi) = (-i)^p \widehat{\varphi}^{(p)}(0) = (-i)^p(fG_0)^{(p)}(0) = (-i)^p\sum_{j=0}^p \binom{p}{j} f^{(j)}(0) G_0^{(p-j)}(0) = 0, \quad p \in \N_0.
$$
Hence, $\varphi \equiv 0$ and, thus, $fG_0 \equiv 0$. Since $G_0$ does not vanish (Lemma \ref{aux-dm}$(i)$), we obtain that $f \equiv 0$.
\end{proof}

In the case of surjectivity, it will frequently come with local extension operators, right inverses for the Borel, respectively, the moment mapping, with a uniform scaling of the constant $h$ determining the Banach spaces under consideration.

\begin{theorem} \label{surj-stiel-modified-dm}
Let $\M$ be a weight sequence satisfying $\sm$,
and let $\A$ be a sequence satisfying $\nq$ and such that $\hA$ is a weight sequence. Then:\par\noindent
(I) Each of the following statements implies the next one:
\begin{itemize}
\item[$(i)$] There exists $a>0$ such that for every $h\ge 1$ there exists a linear and continuous operator $R_h\colon \Lambda_{\M_{+1},h}\to\mathcal{S}^{\hA,ah}_{\M,ah}(0,\infty)$ such that $\mathcal{M}\circ R_h$ is the identity map in $\Lambda_{\M_{+1},h}$.
\item[$(ii)$] There exists $a>0$ such that for every $h\ge 1$ there exists a linear and continuous operator $T_h\colon \Lambda_{\M_{+1},h}\to\mathcal{S}_{\M,ah}(0,\infty)$ such that $\mathcal{M}\circ T_h$ is the identity map in $\Lambda_{\M_{+1},h}$.
\item[$(iii)$] There exists $a>0$ such that for every $h\ge 1$ there exists a linear and continuous operator $U_h\colon \Lambda_{\M_{+1},h}\to C_{\M,ah}(0,\infty)$ such that $\mathcal{M}\circ U_h$ is the identity map in $\Lambda_{\M_{+1},h}$.
\item[$(iv)$] There exists $a>0$ such that for every $h\ge 1$ there exists a linear and continuous operator $V_h\colon \Lambda_{\M_{+1},h}\to\mathcal{A}_{\M_{+1},ah}(\HH)$ such that $\mathcal{B}\circ V_h$ is the identity map in $\Lambda_{\M_{+1},h}$.
\end{itemize}
The following statements are equivalent:
\begin{itemize}
\item[$(i')$] $\mathcal{M}: \mathcal{S}^{\{\hA\}}_{\{\M\}}(0,\infty) \rightarrow \Lambda_{\{\M_{+1}\}}$ is surjective.
\item[$(ii')$] $\mathcal{M}: \mathcal{S}_{\{\M\}}(0,\infty) \rightarrow \Lambda_{\{\M_{+1}\}}$ is surjective.
\item[$(iii')$] $\mathcal{M}: \CM \rightarrow \Lambda_{\{\M_{+1}\}}$ is surjective.
\item[$(iv')$] $\mathcal{B}: \mathcal{A}_{\{\M_{+1}\}}(\HH) \rightarrow \Lambda_{\{\M_{+1}\}}$ is surjective.
\end{itemize}
Any of the statements from (i) to (iv) implies all the conditions from (i') to (iv').\par\noindent
Moreover, the condition
\begin{itemize}
\item[$(v')$] $\ga(\M)>2$.
\end{itemize}
implies (iv).\par\noindent
(II) If $\A$ satisfies in addition the condition $\sm$, then (iv) implies:
\begin{itemize}
\item[$(v)$] There exists $a>0$ such that for every $h\ge 1$ there exists a linear and continuous operator $W_h\colon \Lambda_{\M_{+1},h}\to \mathcal{S}^{(\hA)_{+1},ah}_{\M_{+1},ah}(0,\infty)$ such that $\mathcal{M}\circ W_h$ is the identity map in $\Lambda_{\M_{+1},h}$.
\end{itemize}
\par\noindent
(III) If $\M$ satisfies in addition the condition $\dc$, then we can substitute $\M_{+1}$ by $\M$ in all its appearances, and (iv') implies (v'). So, the six conditions (i')-(v') and (iv) are equivalent.\par\noindent
(IV) If $\M$ and $\A$ satisfy in addition $\dc$, then we can substitute $\M_{+1}$ by $\M$ and $(\hA)_{+1}$ by $\hA$ in all their appearances, and all the ten previous statements (i)-(v), (i')-(v') are equivalent.
\end{theorem}

In the proof of Theorem \ref{surj-stiel-modified-dm} we shall use the following lemma, inspired by \cite{D-E}.
\begin{lemma}(\cite[Lemma 3.6]{DebrouwereetalMoment})\label{inversion-dm}
Let $(c_p)_p \in \C^\N$ and let $G \in C^\infty((-\delta,\delta))$, for some $\delta > 0$, such that $G(0) \neq 0$. Set
$$
b_p = \sum_{j = 0}^p \binom{p}{j} c_j \left ( \frac{1}{G} \right)^{(p-j)} (0), \qquad p \in \N_0.
$$
Then,
$$
\sum_{j=0}^p \binom{p}{j} b_j G^{(p-j)}(0) = c_p, \qquad p\in \N_0.
$$
\end{lemma}

\begin{proof}[Proof of Theorem \ref{surj-stiel-modified-dm}]
$(I)$ $(i) \Rightarrow (ii) \Rightarrow (iii)$: Remark~\ref{rem.ContinuousInclusions} makes these implications obvious.

$(iii) \Rightarrow (iv)$: The map $J_0$ sending every $(c_p)_{p} \in \Lambda_{\M_{+1},h}$ into $((-i)^pc_p)_p$ is a topological isomorphism on $\Lambda_{\M_{+1},h}$, $h>0$. Then it suffices to consider $V_h:=\mathcal{L}\circ U_h \circ J_0$, which, according to Lemma~\ref{inj-lapl-dm}, is linear and continuous from $\Lambda_{\M_{+1},h}$ into $\mathcal{A}_{\M_{+1},bh}(\HH)$ for some $b>0$ independent from $h\ge 1$. Moreover,
$\mu_p(U_h(J_0((c_q)_q))) = (-i)^pc_p$, and so $(V_h((c_q)_q))^{(p)}(0) = i^p \mu_p(U_h(J_0((c_q)_q))) = c_p$ for all $p \in \N_0$, as desired.

$(i') \Rightarrow (ii') \Rightarrow (iii')$ Clear by the inclusions in Remark~\ref{rem.ContinuousInclusions}.

$(iii') \Rightarrow (iv')$: Given $(c_p)_{p} \in \Lambda_{\{\M_{+1}\}}$, pick $\varphi\in C_{\{\M\}}(0,\infty)$ such that $\mu_p(\varphi)=(-i)^pc_p$ for every $p$. Then, $f:=\mathcal{L}(\varphi)\in\mathcal{A}_{\{\M_{+1}\}}(\HH)$ by Lemma~\ref{inj-lapl-dm}, and $f^{(p)}(0) = i^p \mu_p(\varphi) = c_p$ for all $p \in \N_0$.

$(iv') \Rightarrow (i')$: We first claim that $(iv')$ implies that $\ga(\M)>1$ and, so, by Lemma~\ref{lemma.PropertiesGammaIndex}.(iii) we deduce that $(p!)_p\subset\M$. We prove the claim. Observe that, as indicated in Lemma~\ref{extension-dm}, a function $f\in\mathcal{A}_{\{\M_{+1}\}}(\HH)$ can be extended to $\overline{\HH}$ and its restriction to $[-1,1]$ provides a function $f_0$ such that:
\begin{itemize}
\item[(a)] $f_0^{(p)}(0)=f^{(p)}(0)$ for every $p\in\N_0$, and
\item[(b)] $f_0\in\mathcal{E}^{\{\M_{+1}\}}([-1,1])$, the space of functions $\varphi\in C^{\infty}([-1,1])$ such that
    $$
    \sup_{x\in[-1,1],\,p\in\N_0} \frac{|\varphi^{(p)}(x)|}{h^pM_{p+1}}<\infty
    $$
    for suitable $h>0$.
\end{itemize}
Hence, by $(iv')$ the Borel map $\mathcal{B}\colon\mathcal{E}^{\{\M_{+1}\}}([-1,1])\to \Lambda_{\{\M_{+1}\}}$ is also surjective, and a classical result of H.-J. Petzsche~\cite[Thm. 3.5]{Petzsche88} (see also~\cite[Thm. 4.4]{JG-S-S}) proves that $\M_{+1}$ satisfies $(\ga_1)$. By Lemma~\ref{lemma.PropertiesGammaIndex}.(ii), we have $\ga(\M)=\ga(\M_{+1})>1$.

\noindent Consider now the function $G_0$ from Lemma \ref{aux-dm}, and define the linear map $J$ sending every $(c_p)_{p} \in \Lambda_{\M_{+1},h}$, for some $h\ge 1$, into the sequence $(b_p)_p$ given by
$$
b_p = \sum_{j = 0}^p \binom{p}{j} i^jc_j \left ( \frac{1}{G_0} \right)^{(p-j)} (0), \qquad p \in \N_0.
$$
Since this fact will be useful later, we now prove that $J$ is continuous from $\Lambda_{\M_{+1},h}$ into $\Lambda_{\M_{+1},bh}$ for some $b>0$ independent from $h$.
The function $1/G_0$ is holomorphic on a disk with center at 0 and radius larger than 1/2, so there is $C' > 0$ such that $ |(1/G_0)^{(p)}(0)| \leq C' 2^p p!$ for all $p \in \N_0$. Hence,
$$
|b_p| \le  C'|(c_p)_p|_{\M_{+1},h}\sum_{j = 0}^p \binom{p}{j} h^j M_{j+1} 2^{p-j} (p-j)!.
$$
Since $\M$ is $\lc$, we can use~\eqref{eq-conseq-lc-prod}. Also, by the previous argument there exists $B>0$ such that $M_p\ge B^pp!$ for every $p\in\N_0$. So,
\begin{align*}
|b_p| &\le C'|(c_p)_p|_{\M_{+1},h}\sum_{j = 0}^p \binom{p}{j} h^j 2^{p-j} \frac{(p-j)!M_{p+1}}{M_{p-j}}\\
&\le C'|(c_p)_p|_{\M_{+1},h} \left(\left(1+\frac{2}{B}\right)h\right)^pM_{p+1},\ \ p\in\N_0,
\end{align*}
and we are done with $b=1+2/B$.
By assumption, there exists $f\in\mathcal{A}_{\{\M_{+1}\}}(\HH)$ such that $f^{(p)}(0) = b_p$ for all $p \in \N_0$. Lemma~\ref{multiplication-dm} guarantees that $(fG_0)|_{\R} = \widehat{\varphi}$ for some $\varphi \in \mathcal{S}^{\{\hA\}}_{\{\M\}}(0,\infty)$, and Lemma \ref{inversion-dm} implies that
$$
\mu_p(\varphi) = (-i)^p \widehat{\varphi}^{(p)}(0) = (-i)^p(fG_0)^{(p)}(0) = (-i)^p\sum_{j=0}^p \binom{p}{j} b_j G_0^{(p-j)}(0) = c_p, \qquad p \in \N_0,
$$
so we are done.

It is evident that any of the first four statements $(*)$ implies the corresponding $(*')$, and so any of the statements from $(i)$ to $(iv)$ implies all the equivalent conditions from $(i')$ to $(iv')$.

$(v')\Rightarrow(iv)$ Since $\M$ and $\M_{+1}$ clearly share the same gamma index, $(v')$ amounts to $\ga(\M_{+1})>2$. This fact implies in particular (see~\cite[Corollary 3.12 and Remark 3.15]{JimenezSanzSchindlIndex}) that $\M_{+1}\approx\hNN$ for a weight sequence $\NN$ with $\ga(\NN)=\ga(\M_{+1})-1>1$. The condition $\sm$ is stable under passing from $\M$ to $\M_{+1}$, under equivalence, and also under passing from $\hNN$ to $\NN$, so it turns out that $\NN$ satisfies $\sm$ as well, and we can apply~\cite[Theorem 3.3]{JimenezMiguelSanzSchindlSurjectSM} and deduce $(iv)$ for the spaces defined in terms of $\hNN$. Since the equivalence of sequences preserves the spaces of functions or sequences defined by them, this means that $(iv)$ holds as stated.

(II) Observe first that $(iv)$ implies $(iv')$, and as shown in $(iv')\Rightarrow(i')$, we then have $(p!)_p\subset\M$. Since $\A$ satisfies $\sm$, we consider the function $G_0$, and the operator $J$, sending $(c_p)_{p} \in \Lambda_{\M_{+1},h}$, for some $h\ge 1$, into the sequence $(b_p)_p$ as before. By hypothesis, $f:=V_h\circ J((c_p)_p)\in\mathcal{A}_{\{\M_{+1}\}}(\HH)$ is such that $f^{(p)}(0) = b_p$ for all $p \in \N_0$. Now we set $Tf:=(fG_0)|_{\R}$, and define $W_h:=\mathcal{F}^{-1}\circ T\circ V_h\circ J$. According to the behavior described in Lemma~\ref{multiplication-dm} and Proposition~\ref{Fourier-contention-dm} (the hypotheses of the later are easily checked, as $\M_{+1}$ is almost increasing), the map $W_h$ is linear and continuous from $\Lambda_{\M_{+1},h}$ into $\mathcal{S}^{(\hA)_{+1},ah}_{\M_{+1},ah}(0,\infty)$ for some $a>0$ independent from $h$, and $\mathcal{M}\circ W_h$ is the identity map in $\Lambda_{\M_{+1},h}$ by arguing as in $(iv')\Rightarrow(i')$.

(III) Since $\dc$ for the weight sequence $\M$ amounts to $\M_{+1}\approx\M$, the substitution keeps the considered spaces unchanged. Then, it suffices to apply Theorem~\ref{teor.Andreas} to see that $(iv')\Rightarrow(v')$.

(IV) If both  $\M$ and $\A$ satisfy $\dc$, we have $\M_{+1}\approx\M$ and $(\hA)_{+1}\approx\hA$, and so the statements $(i)$ and $(v)$ are equivalent. This fact and the previous implications allow for the conclusion.
\end{proof}

\begin{remark}
The equivalence of the five conditions  $(i')-(v')$ when $\M$ is strongly regular (and so also $\M_{+1}\approx\M$) was already shown in~\cite[Thm. 3.5]{DebrouwereetalMoment}, while the case when $\M$ is $\dc$ is deduced in~\cite[Thm. 6.1.(b) and Thm. 7.2.(b)]{momentsdebrouwere}. The novelty in this situation consists in the equivalence with $(iv)$. One should also note that in~\cite{momentsdebrouwere} the Stieltjes moment problem is also solved for Beurling-like classes, and the existence of global right inverses in both the Roumieu and Beurling classes is characterized. However, the techniques used there seem to heavily depend on the condition $\dc$ (see, for example, Lemmas 3.6.(b) and 3.7.(b) and Proposition 5.1 in~\cite{momentsdebrouwere}), so they are not available under the weaker condition $\sm$. Nevertheless, $\sm$ allows for the construction of the local extension operators $V_{h}$ in $(iv)$, and this has been the motivation for this new insight. Note that, when $\M$ and $\A$ are $\dc$, the construction of local right inverses for $\mathcal{M}$ as the ones in $(i)-(iii)$, with a uniform scaling of the parameter $h$ entering the definition of the corresponding Banach spaces, is new, although $(ii)$ was previously obtained in~\cite{L-S09} when $\M$ is strongly regular.
\end{remark}

\begin{examples}
\item[(i)] As said before, the weight sequence $\M_{q,\sigma}:=(q^{n^\sigma})_{n\in\N_0}$ ($q>1$, $\sigma>1$) satisfies $\dc$ if and only if $\sigma\le 2$; so, former results in the literature about the Stieltjes moment problem in Gelfand-Shilov classes did not cover the case $\sigma>2$. As $\sm$ is satisfied for such $\sigma$, and moreover $\gamma(\M_{q,\sigma})=\infty$, our results prove the moment mapping to be surjective onto $\Lambda_{\{(\M_{q,\sigma})_{+1}\}}$. On the other hand, injectivity does not hold, as Theorem~\ref{inj-stiel-dm} shows.
\item[(ii)] Similarly, the sequence $\M^{\tau,\sigma}=(n^{\tau n^{\sigma}})_{n\in\N_0}$, studied in a series of papers by S. Pilipovi\'c, N. Teofanov and F. Tomi\'c (see~\cite{PilipovicTeofanovTomic} and the references therein), satisfies $\dc$ if and only if $1<\sigma<2$. For $\sigma\ge 2$, surjectivity onto $\Lambda_{\{\M^{\tau,\sigma}_{+1}\}}$ holds again for the Stieltjes moment mapping, because $\sm$ is satisfied and $\gamma(\M^{\tau,\sigma})=\infty$. Again, injectivity fails to hold in this case.
\end{examples}

\section{The situation when $\sm$ fails to hold}\label{sect-Stiel-mom-problem-no-cond}

Our results in the previous section do not apply to rapidly growing sequences not satisfying $\sm$, such as $\M=(q^{p^p})_p$ for $q>1$. In this case we have $\gamma(\M)=\infty$, and since $\sm$ fails, we know by Proposition~\ref{prop.sm.charact.LambdaM+1} that $\mathcal{M}(\CM)\not\subset\Lambda_{\{\M_{+1}\}}$. At the same time, from the estimates~\eqref{eq.bound_Laplace_transform} we see that for every $\varphi\in C_{\M,h}(0,\infty)$, $\zeta\in\HH$ and $p\in\N_0$ we have
$$
|(\mathcal{L}(\varphi))^{(p)}(\zeta)|\le s_{\M,h}^0(\varphi)h^{p+1}M_{p+1}\log\big(\frac{e^2m_{p+1}}{m_p}\big).$$
Taking into account the equalities~\eqref{eq.moments.deriv.Laplace}, it is clear that, if we define $\widetilde{\M}=(\widetilde{M}_p)_p$ as
\begin{equation}\label{eq.defWidetildeM}
\widetilde{M}_p:=M_{p+1}\log\big(\frac{e^2m_{p+1}}{m_p}\big),\quad p\in\N_0,
\end{equation}
we have that the mappings $\mathcal{L}\colon C_{\{\M\}}(0,\infty)\to \mathcal{A}_{\{\widetilde{\M}\}}(\HH)$ and $\mathcal{M}\colon \CM \rightarrow \Lambda_{\{\widetilde{\M}\}}$
are well-defined and continuous (and $\mathcal{L}$ is, as before, injective).

Although a deeper study of this new sequence $\widetilde{\M}$ will be done in Subsection~\ref{subsect.commentsMtilde}, we comment that an $\lc$ sequence $\M$ satisfies $\dc$, respectively $\sm$, if and only if $\widetilde{\M}$ is equivalent to $\M$, resp. to $\M_{+1}$. So, the consideration of $\widetilde{\M}$ makes sense precisely when $\sm$ fails. For example, for $\M=(q^{p^p})_p$ we have that $\widetilde{\M}$ is equivalent to  $(M_{p+1}\cdot p^p)_p$.

\begin{remark}
Suppose $\M$ is a sequence such that $\inf_{p\in\N_0}m_p>0$. Then, as indicated in~\cite[Remark 2.4]{JimenezMiguelSanzSchindlSurjectSM}, the condition $\sm$ can be rephrased as
\begin{equation}\label{eq.sm_modified}
\exists C_1 > 0,\ H > 1\colon \log(m_p)\le C_1H^p,\quad p\in\N_0.
\end{equation}
In particular, this holds if $\M$ satisfies $\lc$; if for such an $\M$ we had $\lim_{p\to\infty}m_p<\infty$, then~\eqref{eq.sm_modified} would be clearly valid. In other words, an $\lc$ sequence not satisfying $\sm$ has to be a weight sequence. In the sequel, we will always consider this situation.
\end{remark}

The following results are straightforward adaptations of Propositions~\ref{Fourier-contention-dm} and~\ref{Fourier-image-dm}, respectively.

\begin{prop}\label{prop.Fourier-contention-sm}
Let $\M$ be a weight sequence such that $\widetilde{\M}$ is almost increasing, and $\A$ be either an almost increasing sequence, or a sequence such that $\liminf_{p\to\infty}A_p^{1/p}>0$ and $\hA$ satisfies $\alg$. Then
\begin{itemize}
\item[$(i)$] $\mathcal{F}\colon\mathcal{S}^{\{\hA\}}_{\{\M\}}(\R)\to \mathcal{S}^{\{\widetilde{\M}\}}_{\{\hA\}}(\R)$ is continuous,
and the same holds at the level of Banach spaces with a uniform scaling of the value $h$.
\item[$(ii)$] The previous statement is valid also for $\mathcal{F}^{-1}$.
\end{itemize}
\end{prop}

\begin{proof}
The proof of Proposition~\ref{Fourier-contention-dm} initially works, except for the last line in~\eqref{eq.estimate-Fourier-middle-intervals}, which cannot be written now. Gathering the valid estimates, we obtain
\begin{equation*}
\sup_{\xi \in \R} |\xi^q\widehat{\varphi}^{(p)}(\xi)|\leq \sum_{j=0}^{\min\{p,q\}} \binom{q}{j}\binom{p}{j} j! 2s_{\M,h}^{\hA,h}(\varphi)h^{p+q+1-2j}(q-j)!A_{q-j}M_{p+1-j}\log\left(\frac{e^2m_{p+1-j}}{m_{p-j}}\right).
\end{equation*}
Since $\widetilde{\M}$ is assumed to be almost increasing, the rest of the proof of Proposition~\ref{Fourier-contention-dm} can be easily adapted.
\end{proof}

\begin{proposition}\label{Fourier-image-without-sm} Let $\M$ be a weight sequence such that $\widetilde{\M}$ is almost increasing, and $\A$ be either an almost increasing sequence, or a sequence such that $\liminf_{p\to\infty}A_p^{1/p}>0$ and $\hA$ satisfies $\alg$. If $\psi \in \mathcal{S}^{\{\widetilde{\M}\}}_{\{\hA\}}(\R)$ and there is $\Psi: \overline{\HH} \rightarrow \C$ satisfying the following conditions:
\begin{itemize}
\item[$(i)$] $\Psi_{|\R} = \psi$.
\item[$(ii)$] $\Psi$ is continuous on $\overline{\HH}$ and analytic on $\HH$.
\item[$(iii)$] $\lim_{\zeta \in \overline{\HH}, \zeta \to \infty} \Psi(\zeta) = 0$,
\end{itemize}
then $\psi \in \mathcal{F}(\mathcal{S}^{\{\hA\}}_{\{\M\}}(0,\infty))$.
\end{proposition}

We also get the analogue of Lemma~\ref{multiplication-dm}.

\begin{lemma}\label{multiplication-without-sm}
Let $\M$ be a weight sequence such that $\widetilde{\M}$ is almost increasing, satisfies $\alg$ and $(p!)_p\subset \widetilde{\M}$. Let $\A$ satisfy $\nq$ and that $\hA$ is a weight sequence, and consider the auxiliary function $G_0$ introduced in Lemma~\ref{aux-dm}.
Then, for any
$f\in\mathcal{A}_{\{\widetilde{\M}\}}(\HH)$ one has $(fG_0)|_{\R} \in \mathcal{F}(\mathcal{S}^{\{\hA\}}_{\{\M\}}(0,\infty))$.
\end{lemma}

\begin{proof}
As in the proof of Lemma~\ref{multiplication-dm}, for $h\ge 1$, $f \in \mathcal{A}_{\widetilde{\M},h}(\HH)$, $x\in\R$ and $p,q\in\N_0$ we get from~\eqref{eq.estimates.deriv.G_0} that
\begin{equation*}
|x^p(fG_0)^{(q)}(x)|
\le C\|f\|_{\widetilde{\M},h}2^p\widehat{A}_p \sum_{k=0}^q\binom{q}{k}h^k\widetilde{M}_{k}2^{q-k}(q-k)!.
\end{equation*}
There exists $B>0$ such that $\widetilde{M}_p\ge B^pp!$ for every $p\in\N_0$. Taking also $\alg$ for $\widetilde{\M}$ into account, we have
\begin{align*}
|x^p(fG_0)^{(q)}(x)|&\le
CC_1^q\|f\|_{\widetilde{\M},h}2^p\widehat{A}_p\widetilde{M}_{q} \sum_{k=0}^q\binom{q}{k}h^k2^{q-k}\frac{(q-k)!}{\widetilde{M}_{q-k}}\\
&\le C\|f\|_{\widetilde{\M},h}\left(2C_1(1+2/B)h\right)^{p+q} \widehat{A}_p\widetilde{M}_{q}.
\end{align*}
Hence,  $(fG_0)|_{\R} \in \mathcal{S}_{\{\hA\}}^{\{\widetilde{\M}\}}(\R)$ and Proposition~\ref{Fourier-image-without-sm} allows us to conclude.
\end{proof}

The next characterization of injectivity is proven in exactly the same way as Theorem~\ref{inj-stiel-dm}, so we omit the details. Observe that if $\widetilde{\M}$ is a weight sequence, Theorem~\ref{inj-borel-dm} guarantees the equivalence of the forthcoming statements $(i)$ and $(ii)$ and, together with the condition $(p!)_p\subset\widetilde{\M}$, we have that the hypotheses of Lemma~\ref{multiplication-without-sm} are satisfied.

\begin{theorem}\label{inj-stiel-without-sm}
Let $\M$ be a weight sequence such that $\widetilde{\M}$ is a weight sequence and $(p!)_p\subset\widetilde{\M}$, and let $\A$ be a sequence satisfying $\nq$ and such that $\hA$ is a weight sequence. Then, the following statements are equivalent:
\begin{itemize}
\item[$(i)$] $\displaystyle \sum_{p = 0}^\infty \frac{1}{\widetilde{m}_p^{1/2}} = \infty$.
\item[$(ii)$] $\mathcal{B}: \mathcal{A}_{\{\widetilde{\M}\}}(\HH) \rightarrow \Lambda_{\{\widetilde{\M}\}}$ is injective.
\item[$(iii)$] $\mathcal{M}: \CM \rightarrow \Lambda_{\{\widetilde{\M}\}}$ is injective.
\item[$(iv)$] $\mathcal{M}: \mathcal{S}_{\{\M\}}(0,\infty) \rightarrow \Lambda_{\{\widetilde{\M}\}}$ is injective.
\item[$(v)$] $\mathcal{M}: \mathcal{S}^{\{\hA\}}_{\{\M\}}(0,\infty) \rightarrow \Lambda_{\{\widetilde{\M}\}}$ is injective.
\end{itemize}
\end{theorem}

\noindent Finally, we present a result of surjectivity for the moment mapping in the absence of $\sm$.

\begin{theorem}\label{surj-stiel-modified-without-sm}
Let $\M$ be a weight sequence and let $\A$ be a sequence satisfying $\nq$ and such that $\hA$ is a weight sequence. With the notation of Theorem~\ref{surj-stiel-modified-dm}.(I), but replacing in every instance $\M_{+1}$ by $\widetilde{\M}$, we have:\par\noindent
(I) The implications $(i)\Rightarrow(ii)\Rightarrow(iii)\Rightarrow(iv)$ and $(i')\Rightarrow(ii')\Rightarrow(iii')\Rightarrow(iv')$ in statement $(I)$ are still valid.\par\noindent
(II) Suppose $\widetilde{\M}$ is a weight sequence, then:\par\noindent
(II.a) $(iv')$ implies $(i')$, and any of the statements $(i)-(iv)$ implies all the equivalent statements $(i')-(iv')$.\par\noindent
(II.b) If $\ga(\widetilde{\M})=\infty$, then $(iv')$ is satisfied.
\end{theorem}

\begin{proof}
The arguments for proving $(I)$ are either evident or rest on the use of the Laplace transform.\par
In $(II.a)$, we need $\widetilde{\M}$ to be a weight sequence in order to apply Petzsche's result, as before, and deduce that $\gamma(\widetilde{\M})>1$. From this point the proof resembles that of the same implication in Theorem~\ref{surj-stiel-modified-dm}.\par
For the proof of statement $(II.b)$, note that, by applying~\cite[Corollary 3.12 and Remark 3.15]{JimenezSanzSchindlIndex}, one knows that $\widetilde{\M}\approx\hNN$ for a weight sequence $\NN$ with $\ga(\NN)=\infty$. We can apply~\cite[Thm. 4.5]{JimenezSanzSchindlSurjectDC} and deduce that $\mathcal{B}: \mathcal{A}_{\{\hNN\}}(\HH) \rightarrow \Lambda_{\{\hNN\}}$ is surjective. Since the equivalence of sequences preserves the spaces of functions or sequences defined by them, we get $(iv')$ as stated.
\end{proof}

\noindent In particular, it can be easily checked that Theorem~\ref{surj-stiel-modified-without-sm}.$(II.b)$ applies for $\M=(q^{p^p})_p$.

\subsection{Some comments on the condition $\sm$ and the sequence $\widetilde{\M}$}\label{subsect.commentsMtilde}

The results in these last paragraphs attempt to shed some light on the new condition $\sm$ and on the properties of the sequence $\widetilde{\M}$, entering our arguments when $\sm$ fails to hold.

Firstly, we are going to construct weight sequences $\M$ satisfying $\sm$ while violating $\dc$, or even not satisfying $\sm$, but with an arbitrarily chosen index $\gamma(\M)\in[0,\infty)$. Indeed, we provide a more general construction which shows that the growth index by V. Thilliez can be prescribed for sequences admitting, in a precise sense, any pre-fixed large growth. This is due to the fact that the gamma index measures not only the rate but also the regularity of growth of the sequence of quotients. Therefore, if $\m$ is chosen to be constant in increasingly larger intervals, we can force its gamma index to vanish even if $\m$ grows significantly between those intervals.

Let us introduce the following notion.

\begin{definition}
A  function $g\colon[0,+\infty)\rightarrow[1,+\infty)$ is called a \emph{growth control function} if it is strictly increasing and satisfies $g(0)=1$ and $\lim_{t\rightarrow+\infty}g(t)=+\infty$.
\end{definition}

\begin{proposition}\label{growththm}
For any given growth control function $g$ and any strictly increasing sequence of integers  $(p_j)_{j\in\mathbb{N}_0}$ with $p_0=0$ there exists a weight sequence $\M=(M_p)_p$ such that
\begin{itemize}
\item[$(i)$] $m_p\le g(p)$ for all $p\in\mathbb{N}_0$,
\item[$(ii)$] $m_{p_j}=g(p_j)$ for all $j\in\N_0$.
\end{itemize}
\end{proposition}

\begin{proof}
We introduce the sequence $\M$ in terms of the corresponding sequence of quotients $\m=(m_p)_p$, so that $M_p=\prod_{i=0}^pm_i$.
Let us now set
$$m_{p_j}:=g(p_j),\;\;\;j\in\mathbb{N}_0,\hspace{15pt}m_p:=m_{p_j}, \;\;\;p_j<p<p_{j+1},\;j\in\mathbb{N}_0.$$
By definition it is immediate that $p\mapsto m_p$ is non-decreasing, $m_0=1$, and $\lim_{p\rightarrow+\infty}m_p=+\infty$. Thus, $\M$ is a weight sequence. Moreover, $(i)$ and $(ii)$ are also clear.
\end{proof}

\begin{proposition}\label{growththmGAMMA}
With the previous notation, if  $\sup_{j\in\mathbb{N}}\frac{p_{j+1}}{p_j}=+\infty$, then  $\gamma(\M)=0$.
\end{proposition}

\begin{proof}
Let $\beta>0$ be arbitrary. If the sequence $\left(\frac{m_p}{(p+1)^{\beta}}\right)_p$ were almost increasing, then
$$\exists\;H\ge 1\;\forall\;0\le p\le q:\;\;\;\frac{m_p}{(p+1)^{\beta}}\le H\frac{m_q}{(q+1)^{\beta}}.$$
When evaluating this condition at $p=p_j$ and $q=p_{j+1}-1\ge p_j$ we obtain $\left(\frac{p_{j+1}}{p_j+1}\right)^{\beta}\le H$. Since $$\frac{p_{j+1}}{p_j+1}\ge\frac{1}{2}\frac{p_{j+1}}{p_j},$$ holds for all $j\ge 1$, in view of our hypothesis we get a contradiction. So, $\gamma(\M)=0$.
\end{proof}

Any sequence of positive real numbers is said to be \emph{dominated by a growth control function $g$} if $(i)$ in Proposition~\ref{growththm} holds, and it \emph{strictly approaches $g$} if in addition $(ii)$ is valid.

\begin{theorem}\label{growththm1}
Given a growth control function $g$, a strictly increasing sequence of integers  $(p_j)_{j\in\mathbb{N}_0}$ with $p_0=0$ and $\sup_{j\in\mathbb{N}}\frac{p_{j+1}}{p_j}=+\infty$, and $\beta\ge 0$, there exists a weight sequence $\M^{\beta}=(M^{\beta}_p)_p$, with corresponding quotients $\m^{\beta}=(m_p^\beta)_p$, such that:
\begin{itemize}
\item[$(i)$] $m^{\beta}_p\le g(p)\cdot(1+p)^{\beta}$ for all $p\in\mathbb{N}_0$,

\item[$(ii)$] $m^{\beta}_{p_j}=g(p_j)\cdot(1+p_j)^{\beta}$ for all $j\in\mathbb{N}_0$,

\item[$(iii)$] $\gamma(\M^{\beta})=\beta$.
\end{itemize}
\end{theorem}

\begin{proof}
Let $\boldsymbol{G}^{\beta}$ denote the Gevrey sequence with index $\beta>0$, and let $\M$ be the weight sequence constructed in Proposition~\ref{growththm} for the function $g$ and the sequence $(p_j)_j$. By Proposition~\ref{growththmGAMMA}, $\ga(\M)=0$. Then, it is straightforward that $\boldsymbol{M}^{\beta}:=\boldsymbol{M}\cdot\boldsymbol{G}^{\beta}$ (termwise product) is a weight sequence and $\gamma(\boldsymbol{M}^{\beta})=\gamma(\boldsymbol{M})+\beta=\beta$.
If we define the growth control function $g\cdot(1+\id)^{\beta}: t\mapsto g(t)\cdot(1+t)^{\beta}$, we immediately see that $\boldsymbol{M}^{\beta}$ strictly approaches $g\cdot(1+\id)^{\beta}$. \end{proof}

\begin{remark}
Suppose the definition of $g$ involves a positive real parameter, say $\ell\in(L,+\infty)$ for some $L>0$, and write $g=g_{\ell}$.
If $g_{\ell}\le g_{\ell_1}$ for any $\ell\le\ell_1$ one obtains an ordered one parameter family of growth control functions $\mathfrak{g}:=\{g_{\ell}: \ell\in(L,+\infty)\}$. Let $\mathfrak{h}:=\{h_{\beta}: \beta\in(B,+\infty)\}$ be a further growth control function family and assume that
\begin{equation}\label{growthvaryrelation}
\forall\;\ell_1>\ell>L\;\forall\;\beta>B\;\exists\;A\ge 1\;\forall\;t\ge 0:\;\;\;g_{\ell}(t)\le(g_{\ell}\cdot h_{\beta})(t)\le Ag_{\ell_1}(t).
\end{equation}
If \eqref{growthvaryrelation} holds, then we say that \emph{$\mathfrak{g}$ uniformly controls $\mathfrak{h}$.}
\end{remark}

Let us apply now this information and Theorem \ref{growththm1} to a concrete situation.

\begin{corollary}\label{growthcounterexample}
For any $\beta\in[0,+\infty)$ there exists a weight sequence $\boldsymbol{M}^{\beta}$ satisfying the following properties:
\begin{itemize}
\item[$(i)$] $\boldsymbol{M}^{\beta}$ satisfies $\sm$.

\item[$(ii)$] $\boldsymbol{M}^{\beta}$ does not satisfy  $\dc$ (and hence $\mg$ either).

\item[$(iii)$] $\gamma(\boldsymbol{M}^{\beta})=\beta$.
\end{itemize}
\end{corollary}

\begin{proof}
Consider the growth control function family $\mathfrak{g}:=\{g_H: H>1, g_H(0):=1, g_H(t):=e^{H^t}, t>0\}$ and we apply Theorem \ref{growththm1} to the growth control function $g_{H_0}$, $H_0>1$ arbitrary but fixed. For $(i)$ note that $\mathfrak{g}$ uniformly controls $\mathfrak{h}:=\{(1+\id)^{\beta}: \beta>0\}$. So, given $H_1>H_0$ there exists $A>0$ such that
$$
m_p^{\beta}=m_p(1+p)^{\beta}\le Ae^{H_1^p},\quad p\in\N_0.
$$
Taking logarithms we easily deduce that the condition~\eqref{eq.sm_modified} holds, and so $\boldsymbol{M}^{\beta}$ satisfies $\sm$.
Concerning $(ii)$ note that $\boldsymbol{M}^{\beta}$ strictly approaches $g_{H_0}\cdot(1+\id)^{\beta}$ and so, in particular, we obtain
$$
\sup_{j\in\mathbb{N}_0}\frac{m^{\beta}_{p_j}}{D^{p_j+1}}
=\sup_{j\in\mathbb{N}_0} \frac{e^{H_0^{p_j}}(1+p_j)^{\beta}}{D^{p_j+1}}=+\infty
$$
for any $D\ge 1$, which proves that $\dc$ fails.
\end{proof}

\begin{corollary}\label{growthcounterexample2}
For any $\beta\in[0,+\infty)$ there exists a weight sequence $\boldsymbol{M}^{\beta}$ satisfying the following properties:
\begin{itemize}
\item[$(i)$] $\boldsymbol{M}^{\beta}$ does not satisfy $\sm$ (and hence $\dc$ either).

\item[$(ii)$] $\gamma(\boldsymbol{M}^{\beta})=\beta$.
\end{itemize}
\end{corollary}

\begin{proof}
We apply Theorem \ref{growththm1} to the growth control function $g_{H}$ defined as
$$
g_H(0):=1,\quad g_H(t):=\exp(e^{H^t}),\ t>0,
$$
for some fixed $H>1$, and obtain $\M^{\beta}$ such that $(ii)$ is satisfied.
Since $\boldsymbol{M}^{\beta}$ strictly approaches $g_{H}\cdot(1+\id)^{\beta}$, in particular we obtain
$$
\sup_{j\in\mathbb{N}_0}\frac{\log(m^{\beta}_{p_j})}{D^{p_j+1}}
=\sup_{j\in\mathbb{N}_0} \frac{e^{H^{p_j}}+\beta\log(1+p_j)}{D^{p_j+1}}=+\infty
$$
for any $D\ge 1$, which proves that $\sm$ fails.
\end{proof}

So, the failure of conditions $\dc$ or $\sm$ has no implication on the value of the index $\ga(\M)$ for general weight sequences, and  the statements of Theorems~\ref{surj-stiel-modified-dm} and~\ref{surj-stiel-modified-without-sm} are meaningful.

We turn now to the study of the sequence $\widetilde{\M}$.

\begin{proposition}\label{prop.study_sm}
Let $\M$ be an $\lc$ sequence.
\begin{itemize}
\item[$(i)$] $\M_{+1}\subset\widetilde{\M}$ holds.
\item[$(ii)$] The following statements are equivalent:
\begin{itemize}
\item[(1)] $\M$ satisfies $\sm$.
\item[(2)] $\M_{+1}$ satisfies $\sm$.
\item[(3)] $\M_{+1}\approx\widetilde{\M}$.
\item[(4)] $\widetilde{\M}$ satisfies $\sm$.
\end{itemize}
\end{itemize}
\end{proposition}

\begin{proof}
$(i)$ Since $m_p\le m_{p+1}$ for every $p\in\N_0$, we have
$\widetilde{M}_p\ge 2M_{p+1}$.\par\noindent
$(ii)$ $(1)\Leftrightarrow(2)$ The condition $\sm$ is clearly preserved by any forward shift of indices. Conversely, if $\M_{+1}$ satisfies $\sm$ the corresponding inequality for $\M$ is obtained by taking into account also the position $p=0$, the only one not covered by the hypothesis.\par\noindent
$(1)\Leftrightarrow(3)$ According to (i), the condition (3) is equivalent to the fact that $\widetilde{\M}\subset\M_{+1}$, which means that there exist $C>0$ and $h>1$ such that for all $p\in\mathbb{N}_0$,
$$\log\left(\frac{e^2m_{p+1}}{m_p}\right)\le Ch^{p},
$$
and this amounts to the condition $\sm$ for $\M$ since
$Ch^{p}-2\le Ch^{p}\le (C+2)h^{p}-2$ for all $p$.\par\noindent
$(1)\Rightarrow(4)$ Note that (1) implies (2) and (3), and recall that $\sm$ is stable under equivalence of sequences.\par\noindent
$(4)\Rightarrow(1)$ In order to give the proof in a compact form we use the following representation for the quotients $m_p=\frac{M_{p+1}}{M_p}$:
\begin{equation}\label{eq.represent_m_p_delta_p}
m_p=\exp\left(\sum_{j=0}^p\delta_j\right),\;\;\;p\in\mathbb{N}_0.
\end{equation}
Thus the values $\delta_p$ are obtained by $\delta_{p+1}=\log(m_{p+1}/m_p)$ for all $p\in\mathbb{N}_0$ and $\delta_0=\log(m_0)$. The fact that $\M$ is log-convex is equivalent to having $\delta_p\ge 0$ for all $p$ and, moreover, recall that
$$\widetilde{m}_p=\frac{\widetilde{M}_{p+1}}{\widetilde{M}_p}= m_{p+1}\frac{\log(e^2m_{p+2}/m_{p+1})}{\log(e^2m_{p+1}/m_p)}= m_{p+1}\frac{2+\delta_{p+2}}{2+\delta_{p+1}},$$
which gives
\begin{equation}\label{smviolatelemmaequ1}
\frac{\widetilde{m}_{p+1}}{\widetilde{m}_p}= \exp(\delta_{p+2}) \frac{(2+\delta_{p+3})(2+\delta_{p+1})}{(2+\delta_{p+2})^2},\;\;\;p\in\mathbb{N}_0.
\end{equation}
By assumption $\widetilde{\M}$ satisfies $\sm$, then
\begin{equation}\label{smviolatelemmaequ3}
\exists\;A\ge 1\;\forall\;p\in\mathbb{N}_0:\;\;\;\frac{\widetilde{m}_{p+1}}{\widetilde{m}_p}\le e^{A^{p+1}}.
\end{equation}
Suppose that $\M$ violates $\sm$, i.e.,
\begin{equation}\label{smviolatelemmaequ2}
\forall\;n\in\mathbb{N}\;\exists\;p_n\in\mathbb{N}:\;\;\; \delta_{p_n+1}=\log\left(\frac{m_{p_n+1}}{m_{p_n}}\right)>n^{p_n+1},
\end{equation}
and we will get a contradiction. By combining \eqref{smviolatelemmaequ1} and \eqref{smviolatelemmaequ3} we obtain the estimate
\begin{align*}
e^{A^{p_n}}&\ge\frac{\widetilde{m}_{p_n}}{\widetilde{m}_{p_n-1}}=\exp(\delta_{p_n+1})\frac{(2+\delta_{p_n+2})(2+\delta_{p_n})}{(2+\delta_{p_n+1})^2}
\ge\exp(\delta_{p_n+1})\frac{4}{(2+\delta_{p_n+1})^2},
\end{align*}
thus $(2+\delta_{p_n+1})^2\ge 4\exp\left(\delta_{p_n+1}-A^{p_n}\right)$ for all $n$. Moreover, note that $\delta_{p_n+1}\ge 1$ for all $n$ by \eqref{smviolatelemmaequ2}, and so $(2+\delta_{p_n+1})^2\le(3\delta_{p_n+1})^2=9\delta_{p_n+1}^2$, hence $$
\frac{9}{4}\ge \exp\left(\delta_{p_n+1}-A^{p_n}-2\log(\delta_{p_n+1})\right),\quad n\in\N.
$$
In addition, \eqref{smviolatelemmaequ2} implies that $\lim_{n\rightarrow\infty}\delta_{p_n+1}=\infty$, consequently, $2\log(\delta_{p_n+1})\le\frac{\delta_{p_n+1}}{2}$ holds for all $n$ sufficiently large. For such $n$, we can write
\begin{equation*}
\frac{9}{4}\ge
\exp\left(\delta_{p_n+1}/2-A^{p_n}\right)>
\exp\left(n^{p_n+1}/2-A^{p_n}\right),
\end{equation*}
where the last estimate follows directly from \eqref{smviolatelemmaequ2}.
Altogether, as $n\rightarrow\infty$ the above estimate yields the  contradiction.
\end{proof}

The next final remarks deal with different growth and regularity properties, or indices, for the sequence $\widetilde{\M}$.

\begin{remark}\label{rem.log_conv_m_p}
We make first some comments on the logarithmic convexity of the sequences under consideration.

Let $\M$ be an arbitrary sequence with the normalization
condition $M_2M_0\ge M_1^2$, in other words $m_0\le m_1$. Then the log-convexity of $\m$
implies the log-convexity of $\M$, since we know that $m_{p+1}/m_p\ge m_p/m_{p-1}$ for all
$p\in\mathbb{N}$, and hence $m_{p+1}/m_p\ge m_1/m_0\ge 1$ for all $p\in\mathbb{N}_0$.
Note that the normalization condition holds, in particular, when $M_0=M_1\le M_2$, and in this situation both $\M$ and $\m$ are also non-decreasing. As we will see in the next remark, the log-convexity of $\m$ can enter some arguments when dealing with properties of $\widetilde{\M}$.  

For an $\lc$ sequence $\M$ whose sequence of quotients is given as in~\eqref{eq.represent_m_p_delta_p}, we have $\delta_p\ge 0$ for all $p$. The log-convexity for $\widetilde{\M}$ amounts then to the fact that the right-hand side in~\eqref{smviolatelemmaequ1} is at least 1. In this case, $\widetilde{\M}$ is a weight sequence provided $\M$ is so,
e.g. via $(i)$ in Proposition~\ref{prop.study_sm} (see the paragraph just before Remark~\ref{rem_propertiesMandstability}).
\end{remark}

\begin{remark}\label{rem.indices_infinity_when_sm_fails}
Observe that, for an $\lc$ sequence $\M$ and according to~\eqref{eq.sm_modified}, the failure of $\sm$ precisely means that
\begin{equation}\label{eq.not_sm_with_lc}
\limsup_{p\to\infty}\left(\log(m_p)\right)^{1/p}=+\infty,
\end{equation}
but the corresponding limit might not exist. As an example, define the strictly increasing sequence $(p_n)_{n\in\N}$ given by
$$
p_1=1;\quad p_{n+1}=1+np_n,\ \ \ n\in\N,
$$
and put $m_0=1$ and
$$
\log(m_p)=n^{p_n}\ \ \ \textrm{whenever }p\in[p_n,p_{n+1}-1].
$$
If $\M$ is the sequence with $M_0=1$ and with sequence of quotients $\m=(m_p)_{p\in\N_0}$, it is clear that $\M$ is a weight sequence and that~\eqref{eq.not_sm_with_lc} is satisfied (by considering the positions $(p_n)_n$), but
$$
\left(\log(m_{p_{n+1}-1})\right)^{1/(p_{n+1}-1)}=n^{1/n}\to 1\ \ \textrm{as }n\to\infty.
$$

For sequences growing in a regular way one usually finds that the indices of O-regular variation for the sequence of quotients $\m$ coincide, see~\cite[Section 3]{JimenezSanzSchindlIndex}. In particular, $\ga(\M)$ (which is the lower Matuszewska index for $\m$) and
$$
\omega(\M)=\liminf_{p\to\infty}\frac{\log(m_p)}{\log(p)}
$$
(the so-called lower order of $\m$) are equal. Suppose that $\sm$ fails and, moreover,
$$
\lim_{p\to\infty}\left(\log(m_p)\right)^{1/p}=\infty.
$$
This easily implies that $\omega(\M)=\infty$. If we suppose also that $\m$ satisfies $\lc$, we have that
\begin{equation*}
\widetilde{m}_p=m_{p+1} \frac{2+\log(m_{p+2}/m_{p+1})}{2+\log(m_{p+1}/m_{p})}\ge m_{p+1}, \quad p\in\N_0,
\end{equation*}
and we can deduce that
$$
\omega(\widetilde{\M})=\liminf_{p\to\infty} \frac{\log(\widetilde{m}_p)}{\log(p)}\ge
\liminf_{p\to\infty}\frac{\log(m_{p+1})}{\log(p)}=\omega(\M)=\infty.
$$
Hence, under mild regularity assumptions, we obtain that also $\ga(\widetilde{\M})=\infty$, and so the item $(II.b)$ of Theorem~\ref{surj-stiel-modified-without-sm} generally proves the surjectivity onto $\Lambda_{\{\widetilde{\M}\}}$ of the moment mapping for such sequences.
\end{remark}

\begin{remark}\label{rem.gamma_index_tildeM}
Let $\M$ be a weight sequence satisfying $\sm$ and such that $\widetilde{M}$ is a weight sequence. According to Proposition~\ref{prop.study_sm}.(ii), and since the gamma index is preserved by forward shifts of indices and by equivalence of weight sequences (see~\cite[Corollary 3.14]{JimenezSanzSchindlIndex}), we deduce that $\gamma(\widetilde{\M})=\gamma(\M_{+1})=\gamma(\M)$. However, we cannot provide any information in case $\sm$ fails to hold.
\end{remark}

\noindent\textbf{Acknowledgements}: The first three authors are partially supported by the Spanish Ministry of Science and Innovation under the project PID2022-139631NB-I00. The research of the fourth author was funded in whole by the Austrian Science Fund (FWF) project 10.55776/PAT9445424.\par

\newpage

\vskip.2cm
\noindent\textbf{Affiliations}:\\
\noindent Javier~Jim\'{e}nez-Garrido:\\
Departamento de Matem\'aticas, Estad{\'\i}stica y Computaci\'on\\
Universidad de Cantabria\\
Avda. de los Castros, s/n, 39005 Santander, Spain\\
Instituto de Investigaci\'on en Matem\'aticas IMUVA, Universidad de Va\-lla\-do\-lid\\
ORCID: 0000-0003-3579-486X\\
E-mail: jesusjavier.jimenez@unican.es\\

\vskip.1cm
\noindent
Ignacio Miguel-Cantero:\\
Departamento de \'Algebra, An\'alisis Matem\'atico, Geometr{\'\i}a y Topolog{\'\i}a\\
Universidad de Va\-lla\-do\-lid\\
Facultad de Ciencias, Paseo de Bel\'en 7, 47011 Valladolid, Spain.\\
ORCID: 0000-0001-5270-0971\\
E-mail: ignacio.miguel@uva.es\\

\vskip.1cm
\noindent Javier~Sanz:\\
Departamento de \'Algebra, An\'alisis Matem\'atico, Geometr{\'\i}a y Topolog{\'\i}a\\
Universidad de Va\-lla\-do\-lid\\
Facultad de Ciencias, Paseo de Bel\'en 7, 47011 Valladolid, Spain.\\
Instituto de Investigaci\'on en Matem\'aticas IMUVA\\
ORCID: 0000-0001-7338-4971\\
E-mail: javier.sanz.gil@uva.es\\

\vskip.1cm
\noindent Gerhard~Schindl:\\
Fakult\"at f\"ur Mathematik, Universit\"at Wien,
Oskar-Morgenstern-Platz~1, A-1090 Wien, Austria.\\
ORCID: 0000-0003-2192-9110\\
E-mail: gerhard.schindl@univie.ac.at
\end{document}